\documentclass[11pt]{article}
\usepackage[utf8]{inputenc}
\usepackage[
papersize={5.5in, 8.5in},
left=0.35in,
right=0.35in,
top=0.35in,
bottom=0.5in]{geometry}
\usepackage{moresize}
\usepackage{graphicx}
\usepackage{amsthm}
\usepackage{rotating}
\usepackage{amsfonts,amssymb,amsmath}
\usepackage[usenames]{color}
\usepackage[T1]{fontenc}

\usepackage{longtable}
\usepackage{array}
\usepackage{float}
\usepackage{booktabs}
\usepackage[dvipsnames,usenames,table,xcdraw]{xcolor}

\newcommand{\gr}{\cellcolor[HTML]{E0E0E0}}

\newcommand\T{\rule{0pt}{2.5ex}}       
\newcommand\B{\rule[-0.7ex]{0pt}{0pt}} 
\newcommand{\comment}[1]{}
\newcolumntype{P}[1]{>{\centering\arraybackslash}p{#1}}

\graphicspath{ {images/} }

\newtheorem{conj}{Conjecture}
\title{\bf 
More certainty in coloring the plane with a forbidden distance interval
}
\author{\bf 
\textcolor[rgb]{0,0.8,0.8}{Jaan Parts} \\
} 
\date{\normalsize \textcolor[rgb]{0,0.8,0.8}{Kazan, Russia, jaan\_parts@.mail.ru}}

\begin{document}

\maketitle

\pagestyle{empty}
\thispagestyle{empty}

\begin{abstract}
In the mysterious and colorful world of chromatic numbers, where there are a lot of unknown, there is an amazing thing. It turns out that for some intervals of forbidden distances on the plane, one can specify the exact value of the chromatic number $\chi$. Two sets of such intervals have been found, for $\chi=7$ and 9. We call them islands of certainty. Here we increase the size of these islands, and add three new ones with $\chi=8$, 12, 13. We also 
conjecture islands for $\chi=14$, 15, 16. Are there islands of certainty for $\chi$=10 or 11? This is still a mystery. Roll up for the Mystery Tour.
\end{abstract}

\section{Introduction}

In the endless ocean of chromatic numbers, it is very difficult to find an island of exact knowledge (except in trivial cases). Usually, instead, we know the bounds on possible values, and improving these bounds is rather slow. Even the sensational breakthrough \cite{grey} of Aubrey de Grey, who narrowed down the set of possible values of the chromatic number of the plane to $\chi\in\{5, 6, 7\}$, can already be considered a page of history\footnote{For more on the history of the coloring problem, see \cite{soi}, and a number of recent developments in this area are outlined in \cite{wes}.}. For other spaces and types of chromatic numbers, the situation is usually even worse.

But even here occasionally there are places where you can feel the earth under your feet. One such place was discovered by Geoffrey Exoo, considering chromatic numbers with an interval of forbidden distances \cite{exoo}.
Namely, for some of these intervals, the exact value of $\chi$ can be found.

\paragraph{Definitions.}
The \textit{chromatic number} of a graph is the minimum number of colors required to color all its vertices so that adjacent vertices take on different colors. Such a coloring is called \textit{proper}. The \textit{chromatic number of the plane} $\chi$ is the chromatic number of an infinite graph, for which the set of vertices coincides with the set of all points of the plane, and the set of edges consists of all pairs of vertices at a unit \textit{distance} from each other. Such a distance is called \textit{forbidden}.

Here we consider $\chi$ as a function of the interval of forbidden distances $[\,1, d\,]$, where $d\ge 1$. By the principle of continuity, $\chi(d)$ is a non-decreasing step function.
As in the case of the traditional chromatic number $\chi(d=1)$, instead of the exact value, only the \textit{bounds} of possible values of $\chi(d>1)$ are usually known. The same is true for the inverse function $d(\chi)$: for each $\chi$, one can specify lower bound $d_{lb}$ and upper bound $d_{ub}$. As a rule, $d_{lb}$ is given by a proper \textit{tiling} of the plane using $\chi$ colors, and $d_{ub}$ can be found using a $\chi$-colorable \textit{finite graph}.

An \textit{island of certainty} is a set 
$d_{min}<d\le d_{max}$ for which the value of $\chi(d)$ is known exactly. The island of certainty arises if $d_{min}<d_{max}$, where $d_{min}(\chi)=d_{ub}(\chi-1)$, and $d_{max}(\chi)=d_{lb}( \chi+1)$.

Note that hereinafter, the letter $d$ denotes the ratio of the upper and lower values of the forbidden distance interval, and not any specific distance. We will also omit subscripts when the meaning is clear from the context.

\paragraph{Background.}

Previously, two islands of certainty were reported. The first such island (with $\chi=7$) was found by Exoo \cite{exoo}. He also proposed
\footnote{
Actually, 
this conjecture was formulated somewhat differently, but in a private correspondence, Exoo confirmed that there is an inaccuracy in the original formulation.}
\begin{conj}[Exoo]
\label{cexoo}
$\chi(d)=7$ for all $d\in(1, \sqrt7/2\,]$.
\end{conj}

Recently, Joanna Chybowska-Sokół, Konstanty Junosza-Szaniawski, and Krzysztof W\k{e}sek\footnote{
For a number of reasons, we refer to this group of authors as W\k{e}sek et al. Firstly, we want to restore some justice: it's a shame to be at the end of the alphabetical list and always hide behind these "et al." (Although we know a good workaround: articles should be written without co-authors.) Secondly, a male surname is usually more stable than a female one. Finally, so shorter.}
discovered  another island of certainty (with $\chi=9$) and put forward \cite{wes} the following 

\begin{conj}[W\k{e}sek et al.]
\label{cwes}
For any integer $k\ge 7$, there exists $d>1$ such that $\chi(d)=k$.
\end{conj}

Or, in our terms, each $\chi\ge 7$ has its own island of certainty.

The purpose of this work is to expand the known islands of certainty and discover new ones. For the most part, we have limited ourselves to studying small $\chi\le 16$. Since the tradition is already visible here, at the end we will propose our conjecture too.

\begin{figure}[!t]
\centering
\includegraphics[scale=0.35]{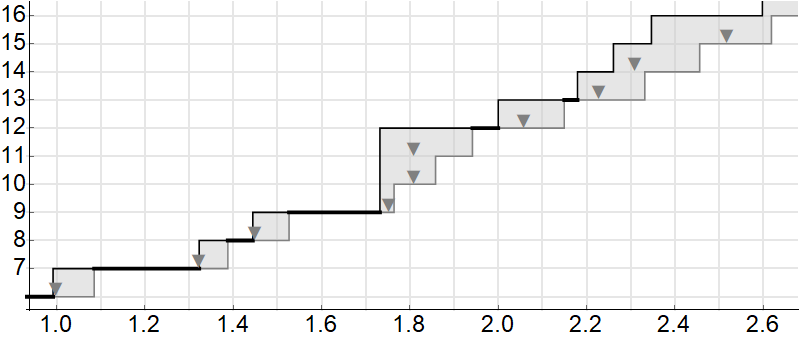}
\caption{Chromatic number $\chi(d)$ for small $d$ values. Black bold horizontal lines mark the islands of certainty. The triangles mark the upper bounds of $d$ obtained by linear extrapolation, as explained in the text.}
\label{fchi}
\end{figure}

\paragraph{Main results.}

Our main results are presented in Fig.~\ref{fchi} and in Table~\ref{tmain}.

For $d<1$ the value $\chi$ is undefined, but can be extended from the side of the lower bounds $d_{lb}$ (obtained by tilings).

\begin{table}[!b]
{
\caption{Main results. Bounds on distance interval $d(\chi)$.
}
\label{tmain}
\smallskip
{
\centering
\footnotesize
\begin{tabular}{@{\;}c|*{7}{>{\!}c<{\!}|}>{\!}c@{\;} } 
\hline
\T\small{$\chi$} & \small{island} & \small{lower}   & \multicolumn{2}{>{\!}c<{\!}|}{\small{island of certainty}} & \small{upper} & \small{clique} & \multicolumn{2}{c}{\small{line}}  \\
\cline{4-5} \cline{8-9}
\B              & \small{status} & \small{bound} & \small{min} & \small{max} & \small{bound} & \small{packing} & \small{pred.} & \small{slope}  \\
\hline
\hline \T
 7 & old  & 0.992076 & 1.085134 & 1.322876 & 1.387777 & 1.414214 & 1.00  & 0.96 \\
 8 & new  & 1.322876 & 1.387777 & 1.444157 & 1.526316 & 1.618034 & 1.323 & 1.72 \\
 9 & old  & 1.444157 & 1.526316 & 1.732051 & 1.764000 & 1.902113 & 1.45  & 1.24 \\
10 & ?   & 1.732051 & 1.764000 & 1.732051 & 1.858032 & 2.000000 & 1.755 & 0.75 \\
11 & ?   & 1.732051 & 1.858032 & 1.732051 & 1.941451 & 2.246979 & 1.81  & 0.8  \\
12 & new  & 1.732051 & 1.941451 & 2.000000 & 2.149463 & 2.569237 & 1.81  & 1.3  \\
13 & new  & 2.000000 & 2.149463 & 2.179449 & 2.331924 & 2.777311 & 2.06  & 1.3  \\
14 & pred & 2.179449 & 2.331924 & 2.260808 & 2.456210 & 2.867455 & 2.23  & 1.3  \\
15 & pred & 2.260808 & 2.456210 & 2.346969 & 2.618615 & 2.909313 & 2.31  & 1.4  \\
\B 
16 & pred & 2.346969 & 2.618615 & 2.598076 &          & 3.151196 & 2.52  & 0.8  \\
\hline
\end{tabular}

}
}
\end{table}

In Table~\ref{tmain} we distinguish four (discovery) states of the island of certainty for each $\chi\in[\,7, 16\,]$: “old” and “new” for confirmed islands (previously known and newly discovered), “pred” and “?” for unconfirmed ones (predicted by linear extrapolation and in doubt). Columns 3 to 7 show the bounds in sequence: $d_{lb}$, $d_{min}$, $d_{max}$, $d_{ub}$, $d_{ub}^*$.

Pairs of estimates $(d_{lb}, d_{max})$ and $(d_{min}, d_{ub})$ repeat each other with a shift by one row, by definition. 
The record estimates of $d_{lb}$ correspond to the tilings shown in Fig.~\ref{fhex} $(k=6, 7, 9, 12, 15)$ and Fig.~\ref{fper} $(k=8, 14, 15)$, with $\chi =k+1$. The record estimates of $d_{ub}$ obtained by graphs (a typical view of which is shown in Fig.~\ref{fgr}) are taken from Table~\ref{tann} ($k=6$) and Table~\ref{texoo} (bottom rows for each $k\neq 9$), with $\chi=k+3$ and $\chi=k$ respectively. The estimates $d_{ub}^{*}$ obtained by point packing are taken from Table~\ref{tcli}, with $\chi=q+3$.
The last two columns of Table~\ref{tmain} show the linear extrapolation parameters obtained from Fig.~\ref{fpred}: the predicted value of $d_{min}$ and the line slope.

The rest of this paper is organized as follows. In Section 2, we consider the basic tools and constructions for obtaining estimates $d(\chi)$. In Sections 3 and 4, we present lower and upper bounds for $d$ obtained from tilings and graphs. In Section 5, we try to predict further progress in refining these bounds. In Section 6, we discuss some difficult cases on which we are hopelessly stuck. Finally, in Section 7 we make some concluding remarks and formulate our conjecture.

\section{Preliminaries}
\paragraph{Poly-chromatic vertices.}

The transition from a single forbidden distance ($d=1$) to a non-zero interval of forbidden distances ($d>1$) greatly simplifies the proof of some facts and changes the numerical estimates. For example, the straight line $\mathbb{R}^1$ can no longer be properly colored in two colors, and a third one is needed. And in the plane $\mathbb{R}^2$, this immediately raises the trivial estimate based on an equilateral triangle from $\chi(d=1)\ge 3$ to $\chi(d>1)\ge 6$.

Indeed, in the case of $d>1$, we can always place the unit triangle on an arbitrarily colored plane in such a way that the $\varepsilon$-neighborhood of one of its vertices will contain at least three colors (will be tri-chromatic), and the $\varepsilon$-neighborhood of another vertex will contain at least two colors (will be bi-chromatic). We call this observation the $(q+3)$-argument, which reflects the increase in $\chi$ compared to the estimate based on a $q$-vertex \textit{clique}. Corresponding but stronger theorems are also available, for example, in Coulson \cite{cou}, Currie-Eggleton \cite{cur}, and W\k{e}sek et al. \cite{wes}. In fact, the latter repeats Coulson's argument: any coloring of the plane in two colors leads to the formation of monochrome stripes or annuli, the length or diameter of which exceeds one.


\paragraph{Graph constructions.}
We will use, with some modifications, the constructions proposed by Exoo \cite{exoo} and by W\k{e}sek et al. \cite{wes}, which we call e-graph and w-graph after their authors, respectively.

By an e-graph we mean a finite graph $r\! \oplus^m\! H$ with $3m^2+3m+1$ vertices located on a hexagonal lattice with step $r=\sqrt{1/a}$ and edges formed by all pairs of vertices at a distance from 1 to $d=\sqrt{b/a}$. Here $a, b\in L=\{u^2+uv+v^2;\; u, v\in\mathbb{Z}_{\ge 0}\}$, the so-called \textit{L\"{o}schian numbers}; $H$ is a 7-vertex \textit{wheel graph} with edges of unit length; and $\oplus^m$ is the Minkowski sum applied to $m$ identical copies of the graph (e.g. $\oplus^3 H=H\oplus H\oplus H$). In our case, the e-graph is bounded by a hexagon with side $r\cdot m$, unlike the original \cite{exoo}, where a rectangle was used. 
It is convenient to omit the normalizing factor $r$ and use the interval of forbidden distances $[\sqrt{a}, \sqrt{b}]$.

In a w-graph, all vertices are located inside the \textit{annulus}, that is, in the area between two concentric circles with radii 1 and $d$. More precisely, $c$ circles are used, each of which has $p$ vertices evenly spaced. Since all $p\cdot c$ vertices are at a forbidden distance from the center, then by placing a tri-chromatic vertex there, we immediately increase the estimate $\chi$ by 3.

Unlike the original constructions of e- and w-graphs, we use both tri- and bi-chromatic vertices, which gives a noticeable improvement in estimates.

\begin{figure}[!t]
\centering
{
\centering
\begin{tabular}{@{}c@{\:}c@{}}
    \includegraphics[scale=0.18]{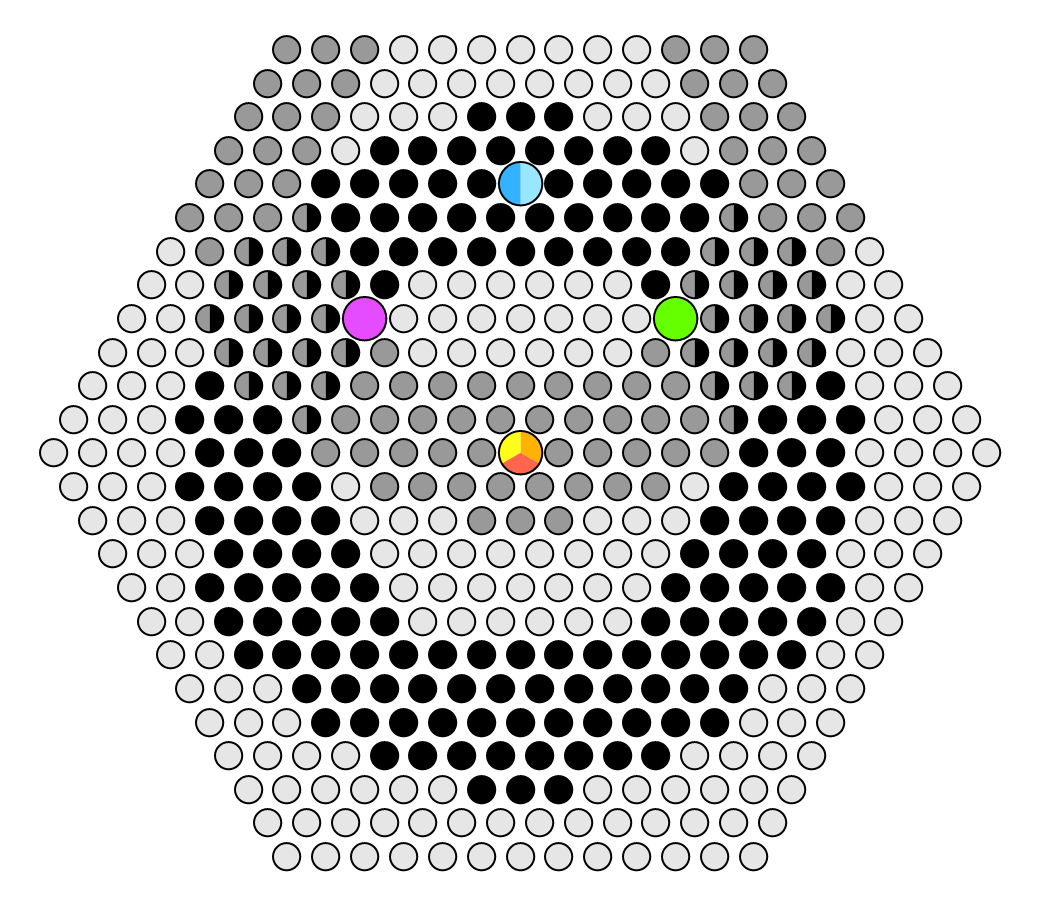} & \includegraphics[scale=0.18]{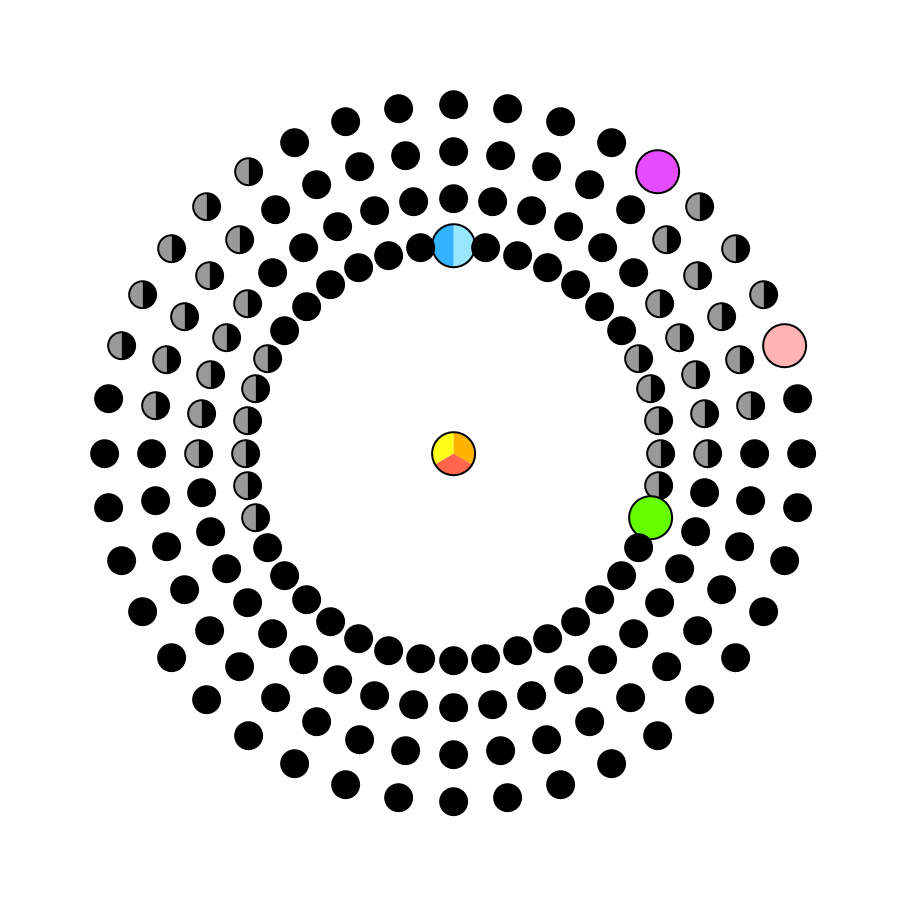} \\
    e-graph, $(m,a,b,q)=(12,27,76,4)$ & 
    w-graph, $(p,c,q)=(40,4,5)$ 
\end{tabular} \par
}
\caption{Examples of e- and w-graphs with $d\approx 1.678$. The vertices of the pre-colored $q$-clique are enlarged. The tri-chromatic vertex is placed in the center of the graph, the bi-chromatic one is above it. Black and dark gray highlight vertices adjacent to tri- and bi-chromatic vertices, respectively.}
\label{fgr}
\end{figure}

Fig.~\ref{fgr} shows examples of e- and w-graphs. The edges are not shown (there are too many of them). Instead, selected vertices adjacent to tri- and bi-chromatic vertices are highlighted. 
The main parameters of the graphs are also listed. Note that the graphs also have other hidden parameters (such as 
the position of the bi-chromatic vertex). 

It can be seen that for the same $d$, the w-graph occupies a smaller area and may have a larger $q$-clique and a smaller number of vertices (especially for small $\chi$). In addition, the $(q + 3)$-argument can significantly reduce the complexity of the computational problem. 


\paragraph{Colors.}
Along with $\chi$, we use the notation $k$ for the number of colors. But one should be careful with this parameter, because in the considered cases, depending on the context, the difference $\chi-k$ can be from 0 to 5. In particular, tiling the plane with $k$ colors gives $d_{lb}$ for $\chi$ of $k$, and tiling the annulus with $k$ colors gives a lower bound on $d_{min}$ for $\chi$ of $k+4$, which can be obtained using a w-graph together with a tri-chromatic vertex.

\paragraph{Tools.}
To calculate $\chi$ of graphs, we use so-called \textit{SAT solvers}. The input information of the solver is a \textit{propositional formula} in CNF format that describes the structure of the graph assuming its proper $k$-coloring (see \cite{heu} for details). To speed up calculations, the vertices of one of the maximum $q$-cliques of the graph are preliminarily colored, which is taken into account in the formula. If the graph is $k$-colorable, then the formula is \textit{satisfiable} ($k$-SAT solution), otherwise it is \textit{unsatisfiable} ($k$-UNSAT solution). The solver may also give no solution in some reasonable time. To determine $\chi$, one need to check several values of $k$.

\section{Tilings}
Our main task is to color all points of the plane with a given number of colors $k$ so as to maximize the interval of forbidden distances $[\,1, d\,]$ for which any two points take on a different color.

In the case of tilings, this problem can be reformulated as follows: for a given $k$, cover the plane with tiles, each of which receives one of $k$ colors and has a \textit{width} (maximum size) not exceeding one, so as to maximize the distance $d$ between the nearest tiles of the same color. This formulation allows us to include in consideration the cases $k<7$ as well.

\paragraph{Identical tiles.}

We start with the case when all tiles have the same shape and are obtained from each other by translation. Such tilings, discussed in detail in \cite{pgray}, demonstrate high efficiency in tiling the plane. Here we briefly repeat some of the results.

\begin{figure}[!t]
\centering
\includegraphics[scale=0.37]{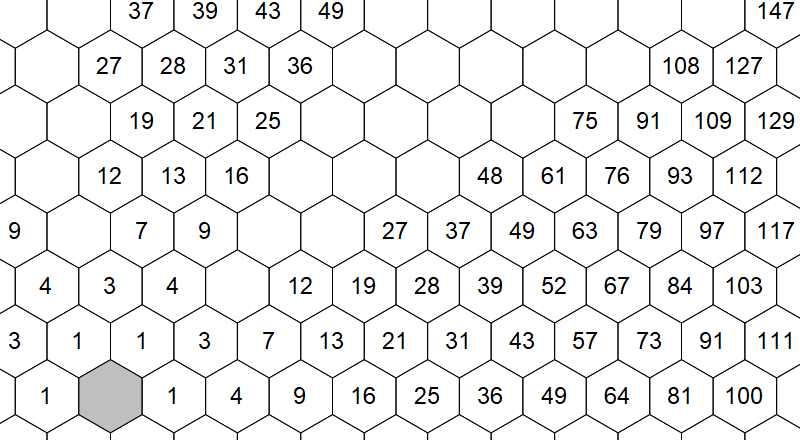}
\caption{Trivial tiling of the plane with regular hexagonal tiles. The base tile is highlighted in gray. Numbered tiles are those closest to the base tile for the given number of colors $k$. Mirrored copies of such tiles are omitted for clarity.}
\label{floe}
\end{figure}

The simplest, and usually the most efficient, are regular hexagons with side length $1/2$. But here tiles of the same color must form a regular hexagonal sublattice, which is possible only if $k\in L$ 
(see Fig.~\ref{floe}). 
In other cases, irregular hexagons give larger value of $d$. Table~\ref{thex} shows the optimal values of $d$ for all $k\in[1,\, 200]$. The cases $k\in L$ are highlighted in gray. It can be seen that an increase in $k$ does not always lead to an increase in $d$. Moreover, for some $k$, a decrease in $d$ is observed. However, this is due to the implicit requirement for tiles of the same color to form a lattice. Relaxing this requirement yields a nondecreasing function $d(k)$, but leaves the question of whether this function is strictly increasing in the case of arbitrary tilings.

\begin{table}[!t]
{
\caption{Minimum distances for lattice-sublattice coloring.}
\label{thex}
\smallskip
\scriptsize

{
\centering
\scriptsize
\begin{tabular}{@{\;}c@{\;\;}|@{\;\;\;}*{9}{>{\!\!}c<{\!\!\!\!}}*{1}{>{\!\!}c}@{\;}}

\hline
\T\B\
\footnotesize{$k$}  & +1 & +2 & +3 & +4 & +5 & +6 & +7 & +8 & +9 & +10 \\
\hline \T
   +0 &\gr 0.00000 & 0.00000 &\gr 0.50000 &\gr 0.86603 & 0.83333 & 0.99208 &\gr 1.32288 & 1.40000 &\gr 1.73205 & 1.65831 \\
  +10 & 1.67416 &\gr 2.00000 &\gr 2.17945 & 2.17945 & 2.18661 &\gr 2.59808 & 2.51979 & 2.47863 &\gr 2.78388 & 2.86282 \\
  +20 &\gr 3.04138 & 3.01579 & 2.94299 & 3.17878 &\gr 3.46410 & 3.39077 &\gr 3.50000 &\gr 3.60555 & 3.42423 & 3.76229 \\
  +30 &\gr 3.90512 & 3.87199 & 3.91379 & 4.00412 & 3.98392 &\gr 4.33013 &\gr 4.27200 & 4.19539 &\gr 4.44410 & 4.41935 \\
  +40 & 4.50035 & 4.65107 &\gr 4.76970 & 4.73559 & 4.60794 & 4.84520 & 4.78005 &\gr 5.00000 &\gr 5.19615 & 5.13667 \\
  +50 & 5.06837 &\gr 5.29150 & 5.12431 & 5.33379 & 5.25115 & 5.53345 &\gr 5.63471 & 5.60196 & 5.47726 & 5.69491 \\
  +60 &\gr 5.76628 & 5.69604 &\gr 5.89491 &\gr 6.06218 & 6.00871 & 5.99567 &\gr 6.14410 & 5.98639 & 6.17729 & 6.20261 \\
  +70 & 6.16402 & 6.41174 &\gr 6.50000 & 6.46932 &\gr 6.50000 &\gr 6.55744 & 6.16188 & 6.54545 &\gr 6.72681 & 6.70343 \\
  +80 &\gr 6.92820 & 6.87989 & 6.81757 &\gr 7.00000 & 6.94304 & 7.02700 & 6.82668 & 6.97285 & 7.15312 & 7.28729 \\
  +90 &\gr 7.36546 & 7.33696 &\gr 7.36546 & 7.40759 & 7.28947 & 7.44444 &\gr 7.56637 & 7.35575 & 7.59056 &\gr 7.79423 \\
 +100 & 7.75030 & 7.70955 &\gr 7.85812 & 7.80525 & 7.88072 & 7.69172 & 7.84502 &\gr 8.00000 &\gr 8.04674 & 8.16091 \\
 +110 &\gr 8.23104 &\gr 8.20462 & 8.10131 & 8.26755 & 8.15518 & 8.30005 &\gr 8.41130 & 8.18005 & 8.33073 & 8.47751 \\
 +120 &\gr 8.66025 & 8.62006 & 8.56481 &\gr 8.71780 & 8.58790 & 8.73717 &\gr 8.76071 & 8.70315 &\gr 8.84590 & 8.70760 \\
 +130 & 8.85303 & 9.03313 &\gr 9.09670 & 9.07217 & 8.97583 & 9.12889 & 9.11716 & 9.11822 &\gr 9.26013 & 9.21474 \\
 +140 & 9.18245 & 9.31939 & 9.35781 &\gr 9.52628 & 9.48928 & 9.45424 &\gr 9.57862 &\gr 9.53939 & 9.59554 & 9.52353 \\
 +150 &\gr 9.65660 & 9.68793 & 9.58933 & 9.71491 & 9.73937 &\gr 9.90430 &\gr 9.96243 & 9.93960 & 9.84962 & 9.99120 \\
 +160 & 9.95752 & 9.99490 &\gr 10.1119 & 9.97235 & 10.1062 & 10.1655 & 10.0539 & 10.2373 &\gr 10.3923 & 10.3581 \\
 +170 &\gr 10.3320 &\gr 10.4403 & 10.3258 & 10.4553 &\gr 10.4762 & 10.4581 & 10.5400 & 10.3539 & 10.4812 & 10.5999 \\
 +180 &\gr 10.6888 & 10.7747 &\gr 10.8282 & 10.8069 & 10.7747 & 10.8542 & 10.3829 & 10.8426 &\gr 10.9659 & 10.9258 \\
\B +190 &10.8944 &\gr 11.0148 &\gr 11.0340 & 10.9440 & 11.1148 &\gr 11.2583 & 11.2265 & 11.1803 &\gr 11.3027 & 11.2188 \\
\hline
\end{tabular}

}
}
\end{table}

\begin{table}[!t]
{
\caption{Classes of coloring for lattice-sublattice scheme.}
\label{tclas}
\smallskip

{
\centering
\footnotesize
\begin{tabular}{@{\;}c@{\;\;}|l@{}}
\hline
\T\B\small{class} & \small{number of colors $k$} \\
\hline
\hline \T
$L+$ & 3, 4, 7, 9, 12, 13, 16, 19, 21, 25, 27, 28, 31, 36, 39, 43, 48, 49, 52, 57, 61, 63, \\
     & 64, 67, 73, 76, 79, 81, 84, 91, 97, 100, 103, 108, 109, 111, 117, 121, 124, 127, \\
     & 129, 133, 139, 144, 147, 151, 156*, 157, 163, 169, 172, 175, 181, 183, 189, \\
     & 192, 193, 196, 199 \\
\hline \T $L-$ & 37, 75, 93, 112, 148, 171 \\
\hline \T
$\overline{L}+$ & 6, 8, 20, 24, 30, 33, 34, 41, 42, 46, 54, 56, 60, 69, 70, 72, 86, 89, 90, 94, 96, \\
     & 99, 105, 110, 114, 116, 120, 126, 131, 132, 136, 142, 143, 149, 152, 154, 155, \\
     & 156*, 160, 162, 166, 168, 174, 177, 180, 182, 186, 195 \\
\hline \T
$\overline{L}-$ & 1, 2, 5, 10, 11, 14, 15, 17, 18, 22, 23, 26, 29, 32, 35, 38, 40, 44, 45, 47, 50, 51, \\
     & 53, 55, 58, 59, 62, 65, 66, 68, 71, 74, 77, 78, 80, 82, 83, 85, 87, 88, 92, 95, 98, \\
     & 101, 102, 104, 106, 107, 113, 115, 118, 119, 122, 123, 125, 128, 130, 134, 135, \\
     & 137, 138, 140, 141, 145, 146, 150, 153, 158, 159, 161, 164, 165, 167, 170, 173, \\
     & 176, 178, 179, 184, 185, 187, 188, 190, 191, 194, 197, 198, 200 \\
\hline
\end{tabular}

}
}
\end{table}


In Table~\ref{tclas}, values of $k$ are divided into four classes $\{L+, L-, \overline{L}+, \overline{L}-\}$ depending on whether they belong to L\"{o}schian numbers 
and whether they lead to an increase (+) in $d$ compared to all previous values\footnote{
The case $k=156$ falls into two classes, $L+$ and $\overline{L}+$, since both variants of tilings exist for it, and both ensure the growth of $d$. And irregular hexagons are more efficient here.}.

\begin{figure}[H]
\centering
{
\centering
\begin{tabular}{@{}c@{\;}cc@{\;}c@{}}
 6 & \includegraphics[scale=0.25]{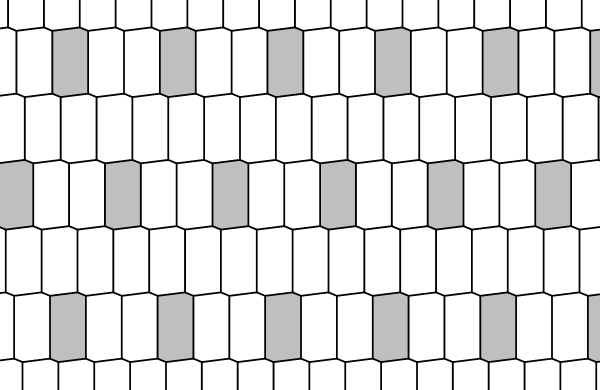}  & 7 & \includegraphics[scale=0.25]{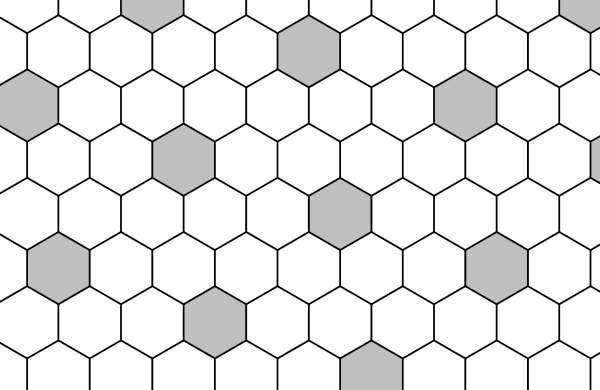}  \\ [2mm]
 8 & \includegraphics[scale=0.25]{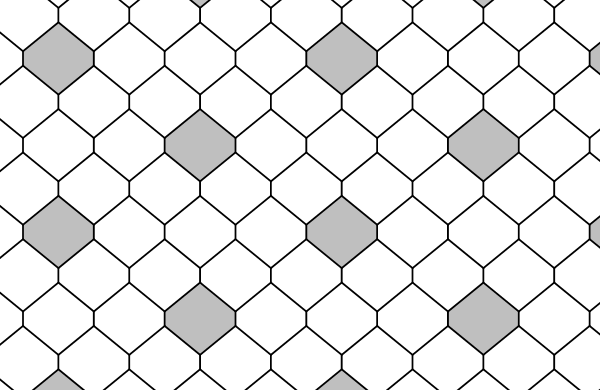}  & 9 & \includegraphics[scale=0.25]{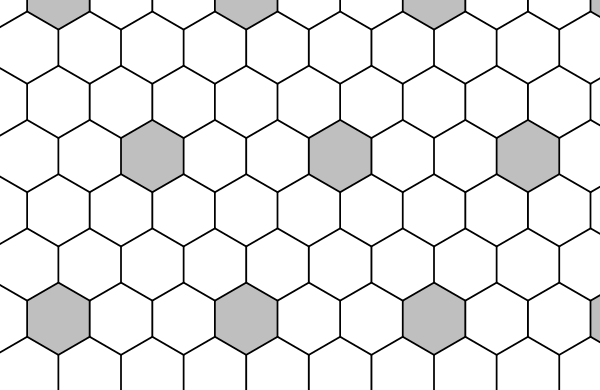}  \\ [2mm]
 10& \includegraphics[scale=0.25]{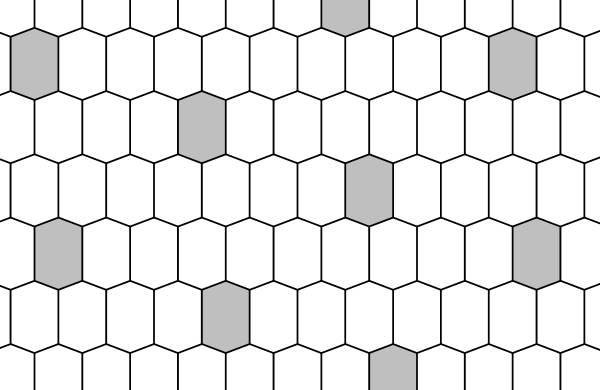} & 11& \includegraphics[scale=0.25]{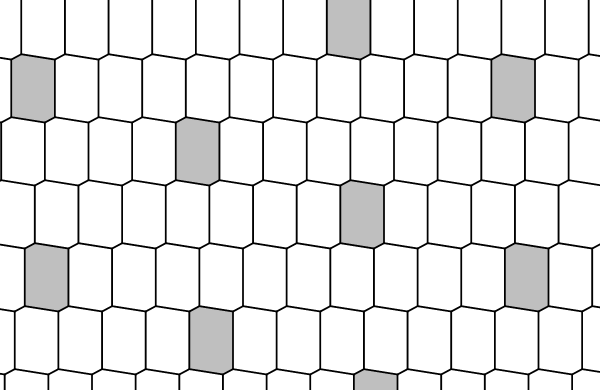} \\ [2mm]
 12& \includegraphics[scale=0.25]{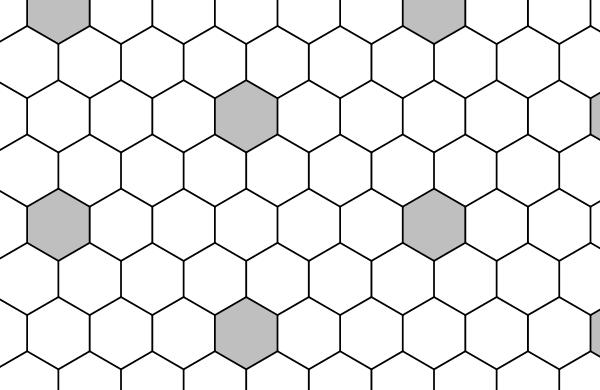} & 13& \includegraphics[scale=0.25]{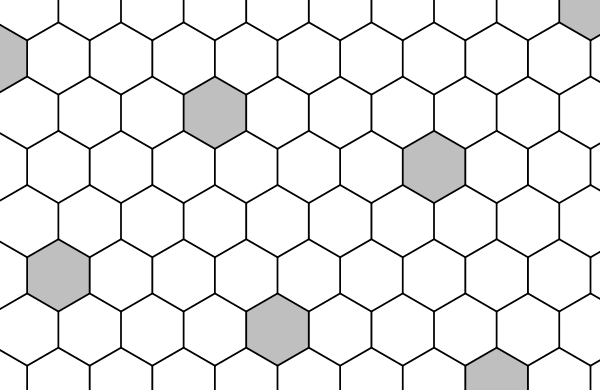} \\ [2mm]
 14& \includegraphics[scale=0.25]{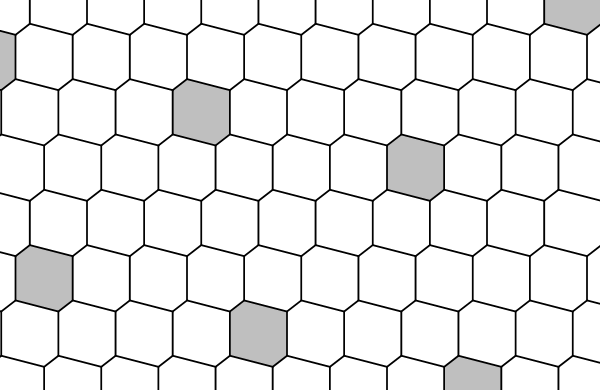} & 15& \includegraphics[scale=0.25]{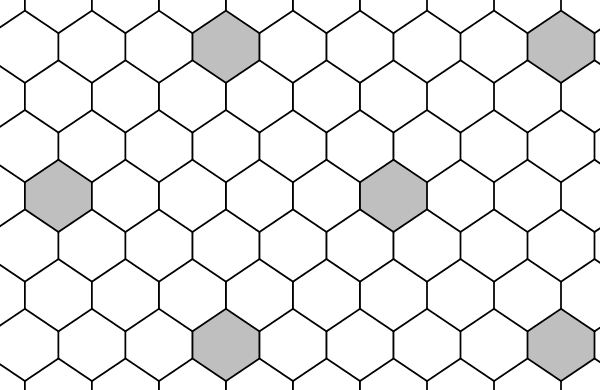} \\
\end{tabular} \par
}
\caption{Optimal hexagonal tilings found in \cite{pgray} for $6\le k\le 15$.}
\label{fhex}
\end{figure}

Fig.~\ref{fhex} shows the best hexagonal tilings identified in [6] for some $k$. Only one of the $k$ colors is shown, the rest are obtained by translation.

\paragraph{Different tiles.}

\begin{figure}[H]
\centering
{
\centering
\begin{tabular}{@{}c@{\;\;}c@{}}
8  & \includegraphics[scale=0.37]{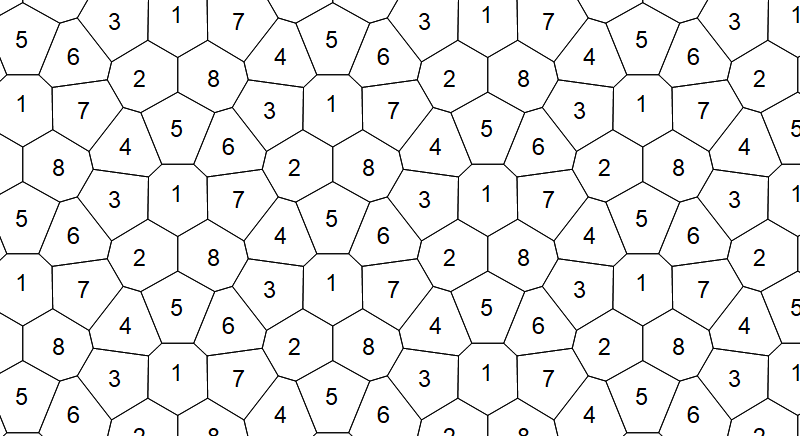}  \\ [4mm]
14 & \includegraphics[scale=0.37]{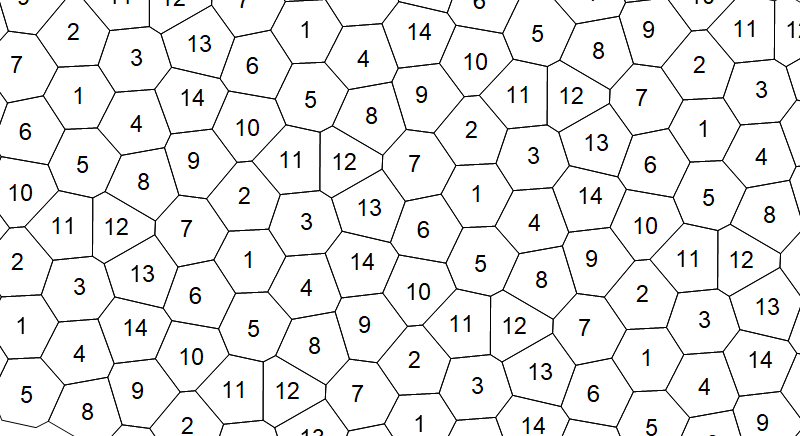} \\ [4mm]
15 & \includegraphics[scale=0.37]{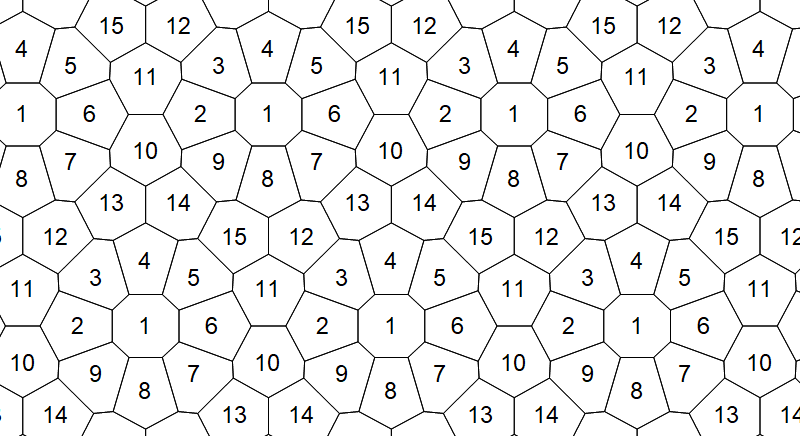} \\
\end{tabular} \par
}
\caption{Non-trivial periodic tilings of the plane with 8, 14, and 15 colors.}
\label{fper}
\end{figure}

We can expand the search for tilings by allowing tiles of different shapes. (To keep the search practical, we however constrain the search by giving each color its own fixed tile shape.) Does this help to increase $d$? Yes, for some $k$ we  found superior tilings, as shown in Fig.~\ref{fper}.

When searching for better tilings, we used the \texttt{HNT} program \cite{sir}, specially developed for solving such problems by Tom Sirgedas. The program allows you to set the tiling area, restrictions on the width and distance between tiles, as well as specify the initial location of the tiles, after which it 
tries to optimize the shape of the tiles. We checked the resulting tilings using the \texttt{FindMaximum} function of the \texttt{Mathematica} software package.

At first glance at the tilings for $k=8$ and 14, one might get the impression that this is a heap of heterogeneous tiles (in the case of $k=15$, the periodicity of the pattern is more obvious). But, if you look closely, you will see the symmetry.

For $k=8$ we got the estimate\footnote{
We give a lot of digits for verification purposes.} $d\approx1.444157346767732$. The tiling includes tiles of three different shapes (considering reflections as equivalent): heptagons $\{1, 2, 8\}$, pentagons $\{3, 5, 7\}$, and hexagons $\{4, 6\}$. If the coloring is not taken into account, the tiling has the symmetry of a regular triangle with respect to the common vertex of the heptagons. All diagonals of pentagons have unit length. Replacing the edge between tiles 2 and 8 with tiles of the ninth color results in a 9-coloring with hexagons.

For $k=14$ we got $d\approx2.260808070967297$. The tiling has third-order rotational symmetry about the center of tiles 1 and 12, and includes tiles of six shapes: $\{1\}$, $\{2, 4, 6\}$, $\{3, 5, 7\}$, $\{8, 11, 13\}$, $\{9, 10 , 14\}$, $\{12\}$.

For $k=15$ we got $d\approx2.346969102448257$. The tiling has horizontal and vertical axes of symmetry passing through the centers of tiles 1 and includes tiles of five shapes: $\{1\}$, $\{2, 4, 6, 8\}$, $\{3, 5, 7, 9\}$, $\{10, 11\}$, $\{12, 13, 14, 15\}$.

None of our attempts to get a more efficient tiling for $k=10$ and 11 were successful (we will return to this problem below).

\paragraph{What's next?}

Is it possible to obtain even more efficient tilings? This question remains open. Note that so far we made some progress whenever we made the tilings more complex: we started with the same shape of all tiles, going from regular hexagons to less symmetrical ones, then we used a different tile shape for each color. One can try to expand the search even more: for example, to abandon the use of a fixed shape of tiles for each color, or to use a different average number of tiles for different colors. However, so far we do not know of any such tiling that would be more efficient in our task.

\newpage
\paragraph{Tiling of annuli}
helps to determine the potential of the w-graph approach. Moreover, the task of annulus coloring can be interesting in itself.

The coloring is constructed in such a way that the inner and outer radii of the annulus numerically coincide with the boundaries of the  forbidden distance interval $[\,1, d\,]$. If the annulus can be tiled using $k$ colors, then any w-graph on this annulus is $k$-colorable. 

W\k{e}sek et al. limited themselves to considering so-called \textit{radial} colorings, in which the annulus is partitioned into sectors bordered by straight lines from its center (and arcs of radii 1 and $d$). It is clear that such a simple construction cannot give estimates $d>2$, that is, it is not applicable for large $k$.

\begin{figure}[H]
\centering
{
\centering
\includegraphics[scale=0.2]{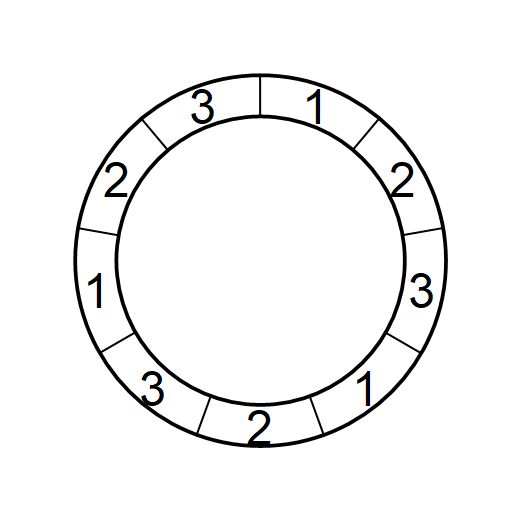}\;\;\;\includegraphics[scale=0.2]{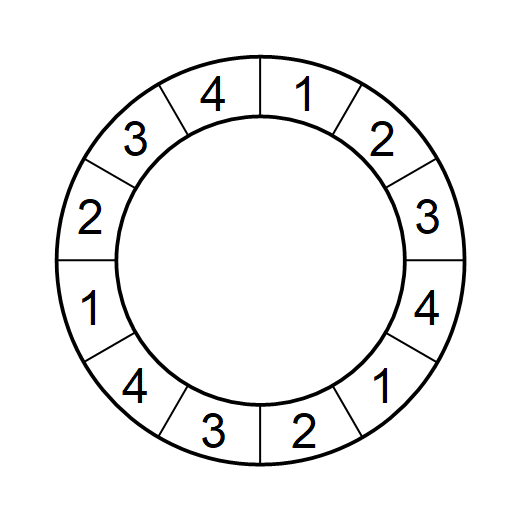}\;\;\;\;\;\,\includegraphics[scale=0.2]{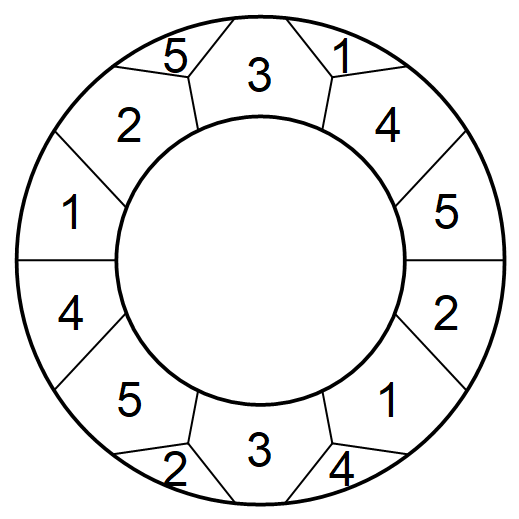}\;\;\, \\[2pt]
\includegraphics[scale=0.2]{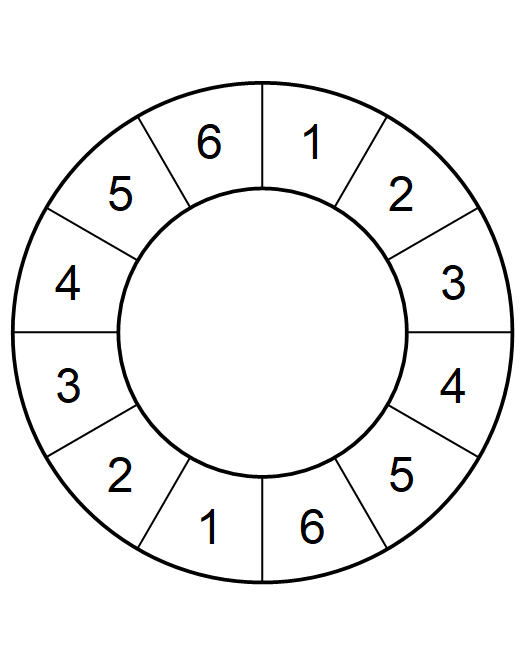}\,\includegraphics[scale=0.2]{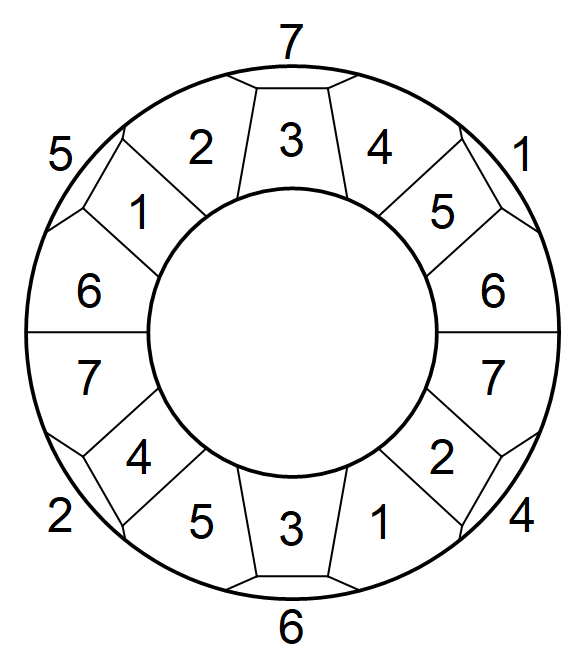}\,\includegraphics[scale=0.2]{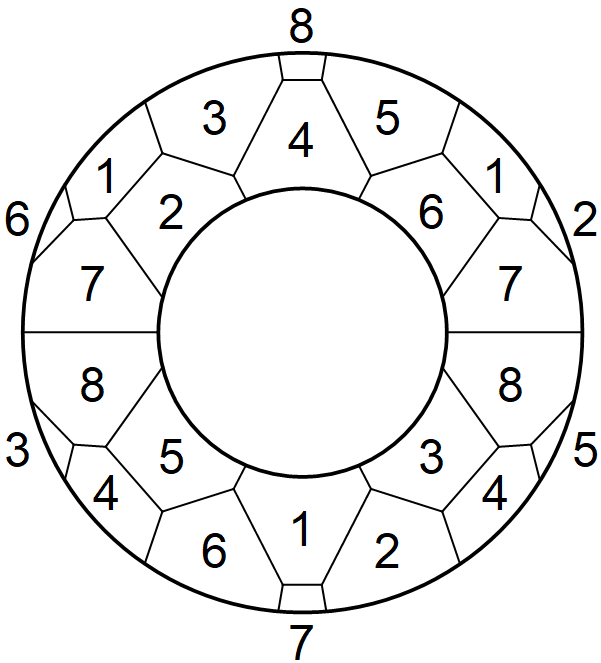} \\[0pt]
\raggedright
\includegraphics[scale=0.2]{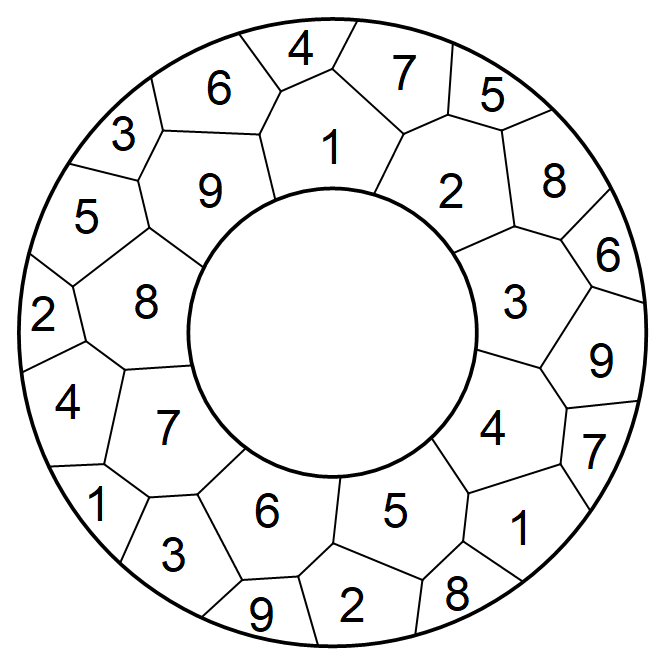}\;\;\; \includegraphics[scale=0.2]{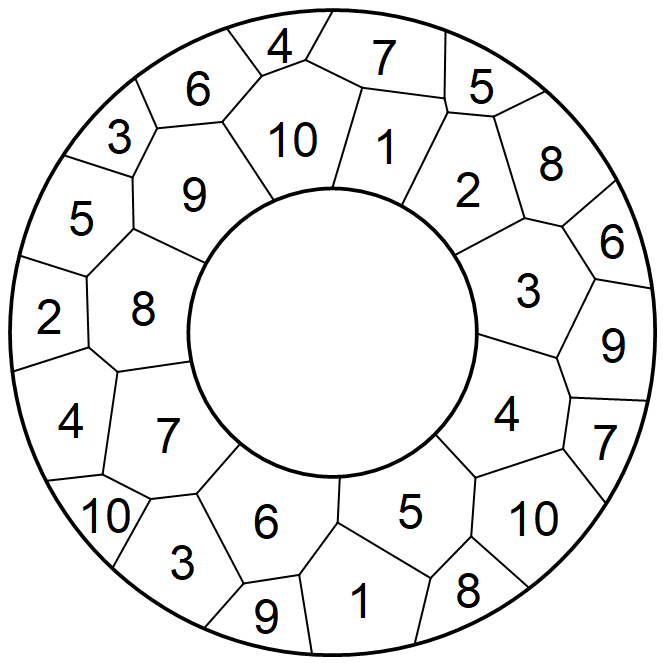} \\[0pt]
\raggedleft
\includegraphics[scale=0.2]{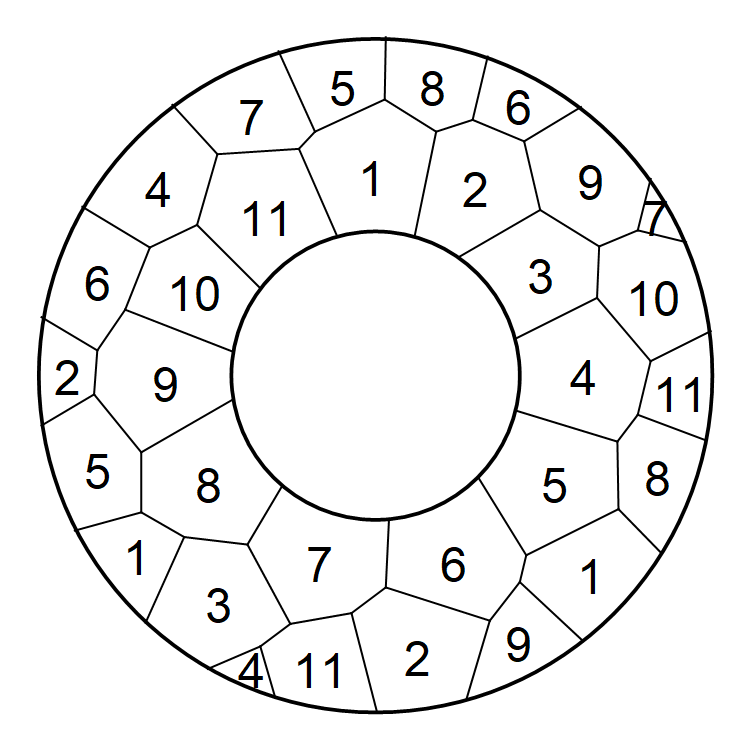}\;\; \includegraphics[scale=0.2]{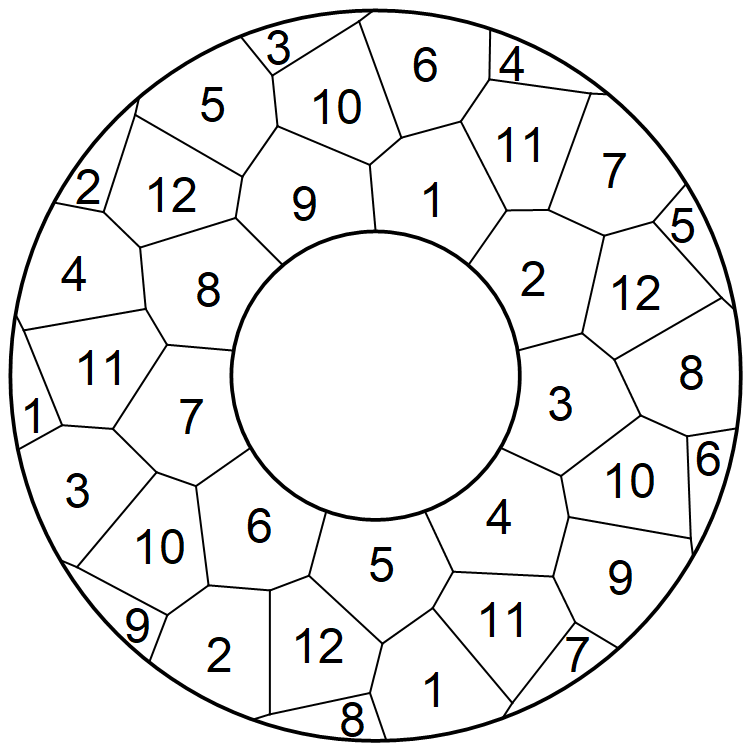}\; \\[0pt]
}
\caption{Coloring of annuli for $3\le k\le 12$.}
\label{fann}
\end{figure}

Fig.~\ref{fann} shows the tilings of the annuli with the largest values of $d$ that we were able to obtain. All annuli are drawn to the same scale and have the same inner radius (although optical illusion convinces us otherwise).

As in the case of the plane, to find the optimal annulus tiling, we used the \texttt{HNT} program and the optimization function of the \texttt{Mathematica} software package.
The refined values of $d$ are shown in the column "arbitrary" of Table~\ref{tann}. For comparison, the "radial" column shows the estimates obtained using the radial coloring.

A radial coloring turned out to be the best for $k=3, 4$, and 6. In other cases, a tiling with multiple "floors" (distances from the center) gives better results. In most cases, we have obtained symmetrical colorings. For $k=10$ and 11, we got an asymmetric coloring by taking the initially symmetrical tiling with $k=9$ and $d\approx 2.175091$, and manually adding tiles of new colors. For these two tilings, the values of $d$ found in the \texttt{HNT} program were not checked (marked with an asterisk in the Table~\ref{tann}).

As tilings show, compared to e-graphs, the w-graph approach has fundamental limitations. For example, for $k=3$ (or $\chi=7$) it is impossible to get $d<2\sin(2\pi/9)\approx 1.285575$, while with the help of e-graphs this barrier can be easily overcome. However, our estimates do not unambiguously indicate the inefficiency of W\k{e}sek's approach for larger $k$. 

We expect that some our estimates of $d$ can be improved.
It is noteworthy that for $k=6$ we failed to improve the radial coloring with $d=\sqrt3$. As you see, 
this value of $d$ is a tough nut to crack. The same difficulties arise when tiling the plane.

\section{Graphs}
\paragraph{Point packing.}

The simplest estimate of the upper bound of $d(\chi)$ can be obtained 
from a set of $q$ points such that the distance between any two points falls within the forbidden interval $[\,1,d\,]$, giving a $q$-clique.
Thus, the problem is reduced to finding the most dense packing of points. 

Although such an estimate has almost no practical use, the point packing problem is interesting in itself. It is similar to the ball-packing problem, but has its own peculiarities: the width of the finite region covering the centers of the balls (and not the entire balls) is minimized.

\begin{figure}[H]
\centering
{
\centering
\includegraphics[scale=0.28]{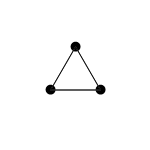}\includegraphics[scale=0.28]{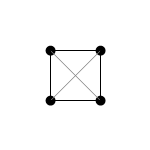}\includegraphics[scale=0.28]{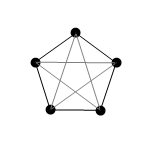}\includegraphics[scale=0.28]{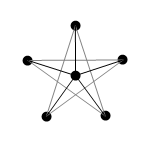}\includegraphics[scale=0.28]{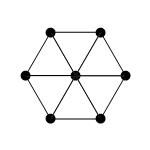}\includegraphics[scale=0.28]{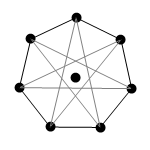} \includegraphics[scale=0.28]{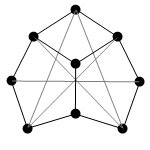} \includegraphics[scale=0.28]{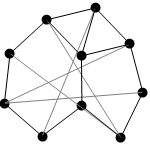} \\[2pt]
\includegraphics[scale=0.28]{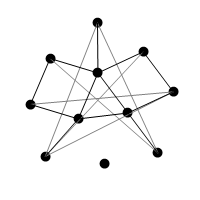}\includegraphics[scale=0.28]{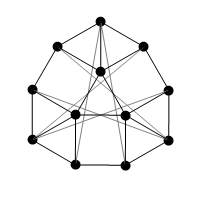}\includegraphics[scale=0.28]{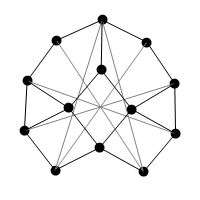}\includegraphics[scale=0.28]{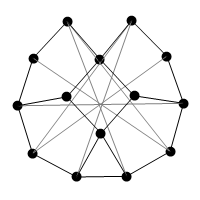} \includegraphics[scale=0.28]{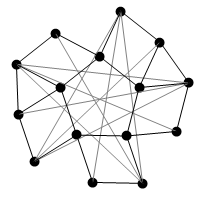} \includegraphics[scale=0.28]{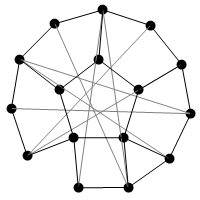}
\includegraphics[scale=0.28]{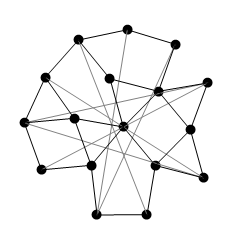}\includegraphics[scale=0.28]{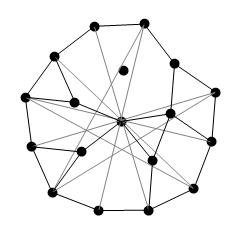}\includegraphics[scale=0.28]{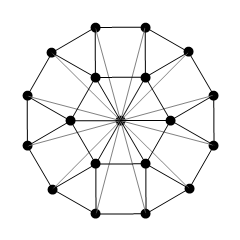}\includegraphics[scale=0.28]{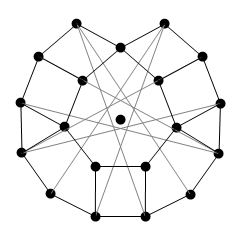}\includegraphics[scale=0.28]{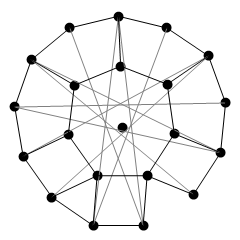}
\includegraphics[scale=0.28]{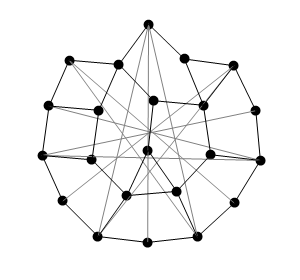}\;\includegraphics[scale=0.28]{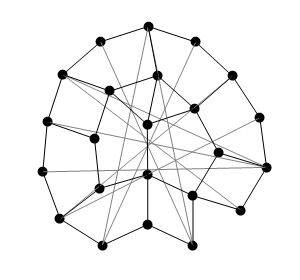}\;\includegraphics[scale=0.28]{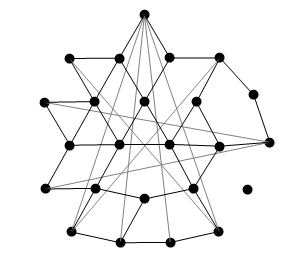}\;\includegraphics[scale=0.28]{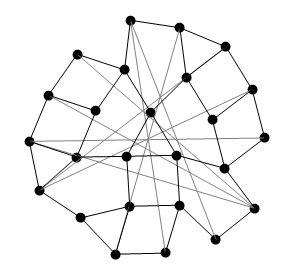}
\includegraphics[scale=0.28]{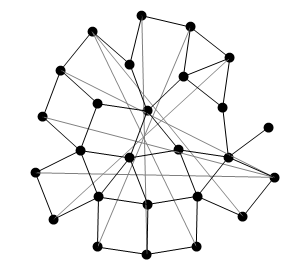}\;\includegraphics[scale=0.28]{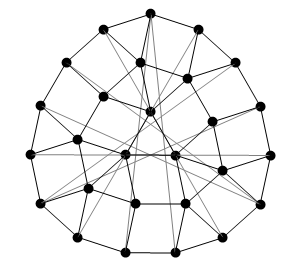}\;\includegraphics[scale=0.28]{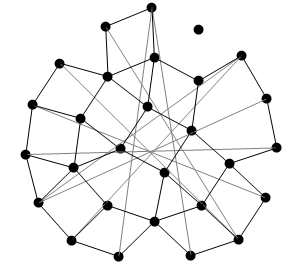}\;\includegraphics[scale=0.28]{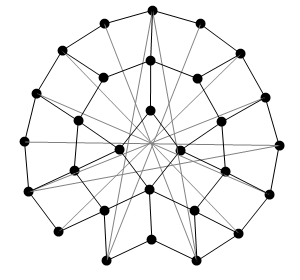}
\includegraphics[scale=0.28]{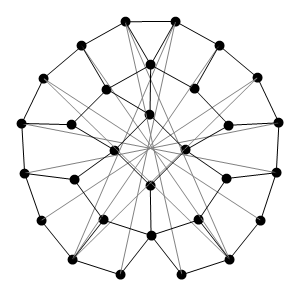}\;\includegraphics[scale=0.28]{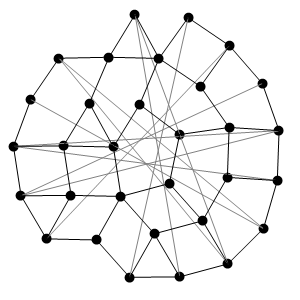}\;\includegraphics[scale=0.28]{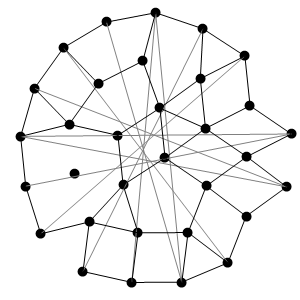}\;\includegraphics[scale=0.28]{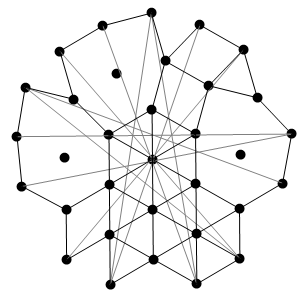} \\[3pt]
\includegraphics[scale=0.28]{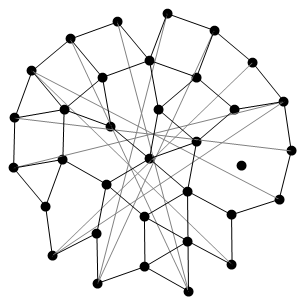}\;\includegraphics[scale=0.28]{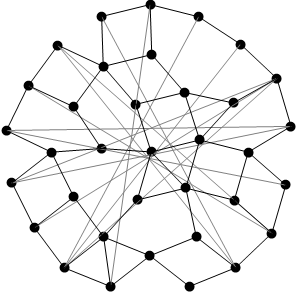}\;\includegraphics[scale=0.28]{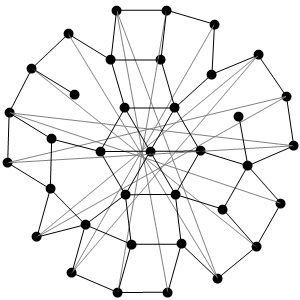}\;\includegraphics[scale=0.28]{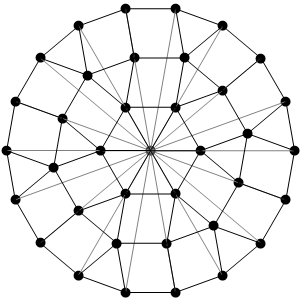}
}
\vspace{-4mm}
\caption{Packing of $q$-cliques for $3\le q\le 37$. The minimal 
and maximal 
distances between vertices are shown in black and gray, respectively.}
\label{fcli}
\end{figure}

The optimization results for some small $q$ are shown in Fig.~\ref{fcli} and Table~\ref{tcli}. It is noteworthy that, as a rule, the point configurations that first come to mind are far from optimal. We used a simple algorithm, but even that proved to be much more successful than trying to find the best solution by hand (with the exception of some symmetrical cases).

\begin{table}[!b]
{
\caption{Minimal width $d$ of $q$-clique for $q\le 40$.
}
\label{tcli}
\smallskip
\scriptsize
{
\centering
\begin{tabular}{@{\;}c@{\;\;}|@{\;\;\;}*{9}{>{\!\!}c<{\!\!\!\!}}*{1}{>{\!\!}c}@{\;}}
\hline
\T\B\
\footnotesize{$q$}  & +1 & +2 & +3 & +4 & +5 & +6 & +7 & +8 & +9 & +10 \\
\hline \T
    +0 & 0.00000 & 1.00000 & 1.00000 & 1.41421 & 1.61803 & 1.90211 & 2.00000 & 2.24698 & 2.56924 & 2.77731 \\
   +10 & 2.86745 & 2.90931 & 3.15120 & 3.31583 & 3.46965 & 3.59016 & 3.75373 & 3.83735 & 3.86370 & 4.09079 \\ 
   +20 & 4.21604 & 4.36509 & 4.47656 & 4.57823 & 4.70821 & 4.77604 & 4.80451 & 5.02067 & 5.09162 & 5.18583 \\
\B +30 & 5.32165 & 5.41447 & 5.49343 & 5.56799 & 5.69445 & 5.73008 & 5.75877 & 5.96806 & 6.02541 & 6.15823 \\
\hline
\end{tabular}

}
}
\end{table}

The calculations were carried out using the \texttt{Mathematica} program. The algorithm is based on the optimization procedure, where the variables are the coordinates of the points, the parameter to be minimized is the largest distance between them, and the constraints are the minimum distance (greater than or equal to one).

The optimization procedure was repeated many times with different initial values of the coordinates, not necessarily satisfying the constraints. At the first stage, we took random initial coordinates in some area. At the next stage, we made small changes to the coordinates that corresponded to the best current solution. Two strategies for changing coordinates were used alternately. The first strategy was to introduce a large ($\pm 2$) random correction to the coordinates of one or several randomly selected vertices. With the second strategy, the coordinates of all vertices were changed by a small random value ($\pm 0.05$). In each case, several hundred iterations were performed. The essence of applying different strategies is to simultaneously solve the problems of approaching the optimum and jumping out of local minima.

The considered algorithm reveals weak convergence, which can be explained by the presence of many local minima in the search area; only a small fraction of iterations leads to an improvement in the result. 
However, the algorithm is quite efficient. The results show that i) as $q$ increases, the symmetry of the packings is usually broken; ii) a fairly common pattern is concentric circles of points; iii) in contrast to tilings, there is a monotonic increase in the estimate of $d$ with the number of colors $q$.

\begin{figure}[!t]
\centering
\includegraphics[scale=0.35]{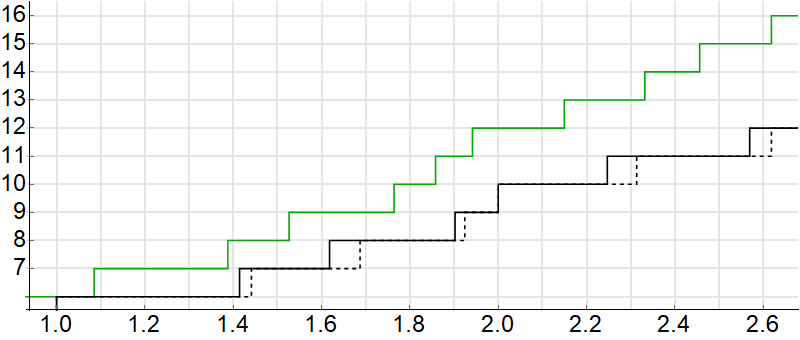}
\caption{The lower bounds on $\chi(d)$ obtained by packing the points of the $q$-clique (lower black curve) and by e- and w-graphs (upper green curve). The dotted line shows an estimate based on a $q$-clique on a hexagonal lattice with step $1/\sqrt{91}$, which turns out to be efficient for many e-graphs.}
\label{flow}
\end{figure}

Fig.~\ref{flow} shows lower bounds on $\chi$ based on point packing given the $(q+3)$-argument. For comparison, estimates based on e- and w-graphs are also shown: it can be seen that with increasing $d$ they diverge more and more.

\paragraph{W\k{e}sek graphs.}

We tried to play with the w-graph parameters, resulting in a slight improvement in numerical estimates compared to the results \cite{wes} of W\k{e}sek et al. (see Table~\ref{tann}, where the latter are listed in the “previous best” column). In the original, W\k{e}sek et al. gave "exact" values of the bounds using trigonometric expressions, but this does not make much sense, since there are no fundamental difficulties in further improving these results. 
We limited ourselves to approximate values. When searching by $d$, we used a step of 0.001 or 0.01 depending on $k$.

We noticed that the dependence of the estimate $d$ on the number of circles $c$ and the number of vertices on the circle $p$ is extremely non-monotonic with narrow local minima. For example, for $k=4$ and $d=1.457$, the local optimum is observed at the following values of $p$: 716, 824, 932, 974, 1040, $\dots$ (we don't see patterns here). For intermediate values of $p$, the estimate of $d$ increases.

A uniform arrangement of the radii of the circles (inside the annulus) is also generally not optimal. For example, for $k=6$ with three circles and $d=1.786$ the optimal radius of the middle circle is about 1.54 ($p=953, 1007,\dots$). With a uniform arrangement of circles, some of them can be discarded. So, for the graph with parameters $(k, d, p, c, q)=(6, 1.764, 294, 12, 4)$, it is possible to remove 6 out of 12 circles without loss of $k$, with numbers 2, 3, 4, 5, 9, 11 (if we count from the center).

\begin{table}[!t]
{
\caption{Estimates of the distance interval $d$ based on coloring of annuli. 
}
\label{tann}
\smallskip
\smallskip
{
\centering
\footnotesize
\begin{tabular}{
@{\;}c@{\;}|
*{3}{>{\!\!}c<{\!\!}|}*{1}{>{\!\!\!}c<{\!\!\!}|}*{1}{>{\!}c<{\!}|}*{1}{>{\!\!\!}c<{\!\!\!}|}*{2}{>{\!\!}c<{\!\!}|}
*{1}{>{\!}c<{\!}|}@{\;}r@{\;}|
*{1}{>{\!}c<{\!}|}>{\!}c@{\;}
}
\hline
\T
\small{$k$}   & \multicolumn{2}{c|}{\small{annulus tiling}}   & \multicolumn{3}{@{}c@{}|}{\small{previous best}} & \multicolumn{5}{@{}c@{}|}{\small{our best w-graph}} & \multicolumn{2}{c}{\small{line}} \\
\cline{2-13}
\B
& \small{radial} & \small{\!arbitrary\!
} & \small{$d$} & \small{$p$} & \small{$c$} & \small{$d$} & \small{$p$} & \small{$c$} & \small{$q$} & \small{time, s} & \small{pred.} & \small{slope} \\
\hline
\hline \T
  3 & 1.285575 & 1.285575 & 1.28599 & 1300 & 2 & 1.2856&   18 &  1 & 2 & 0      &  &  \\
  4 & 1.414214 & 1.414214 & 1.47145 &  180 & 2 & 1.457 &  716 &  2 & 3 & 1.9    & 1.452 & 0.52 \\
  5 & 1.618034 & 1.691392 & 1.71433 &  180 & 3 & 1.696 &  388 &  3 & 4 & 2.2    & 1.685 & 0.54 \\
  6 & 1.732051 & 1.732051 & 1.82843 &  120 & 3 & 1.764 &  294 & 12 & 4 & 23     & 1.749 & 0.72 \\
  7 & 1.801938 & 1.847759 & 2.01176 &  120 & 3 & 1.893 &  363 &  4 & 4 & 182    & 1.869 & 1.23 \\
  8 & 1.847759 & 1.940393 &         &      &   & 1.975 &  232 & 12 & 5 & 15813  & 1.927 & 1.91 \\
  9 & 1.879385 & 2.175091 &         &      &   & 2.224 &  296 & 12 & 6 & 235    & 2.200 & 1.17 \\
 10 & 1.902113 & 2.23649* &         &      &   & 2.380 &  259 & 12 & 7 & 1835   & 2.357 & 1.27 \\
 11 & 1.918986 & 2.33483* &         &      &   & 2.560 &  136 &  6 & 7 & 137504 & 2.498 & 1.18 \\
 12 & 1.931852 & 2.532089 &         &      &   & 2.687 &  150 &  6 & 8 & 139626 & 2.615 & 1.71 \\
\hline
\end{tabular}

}
}
\end{table}


Initially, we tried to apply a brute force method, reducing $d$ only due to a significant increase in $p$ and $c$. But, as practice has shown, an exhaustive search by these parameters turned out to be more effective. Of course, in this case, it was necessary to check a much larger number of graphs, but the total time for checking them still turned out to be much less compared to what was spent on large graphs in the brute force method.
Table~\ref{tann} shows the time of checking single graph with the found optimal parameters using the SAT solver \texttt{glucose}.

\paragraph{Exoo graphs.}

In \cite{exoo}, Exoo used the pair $(a, b)=(25, 43)$ on a 203-vertex graph with $k=6$, which roughly corresponds to our $\oplus^8 H$ graph and gives a $k$-UNSAT solution on $d=\sqrt {43/25}\approx 1.311488$. Checking it with the \texttt{glucose} solver takes about 1000 seconds. Using the pair $(36, 52)$ on $\oplus^8 H$, we get $d=\sqrt{52/36}\approx 1.201850$, but this takes about 
500\,000 seconds, or almost a week. But if we use a modified 
graph with tri- and bi-chromatic vertices, the verification time is reduced to about 10 seconds.

We continued the hunt for records, successively achieving lower and lower estimates of $d$ for several values of $k$.
Simplistically, we adhered to the following scheme for finding the minimum values of $d=\sqrt{b/a}$. We start with a pair of small initial values $(a, b)$, forcing $k$-UNSAT on some suitable graph $\oplus^m H$. With a fixed $a$, the value of $b$ is gradually decreased until a $k$-SAT solution is reached (thus, the preceding $k$-UNSAT solution gives the record value of $d$). Further, the resulting estimate $d$ acts as a threshold, and at approximately constant $d$, the values of $a$ and $b$ are gradually increased until a $k$-UNSAT solution is found. Then the value of $a$ is fixed again and the cycle repeats.

The results are shown in Table~\ref{texoo}. Here, for each $k$, the parameters of several graphs that give a $k$-UNSAT solution are shown, including the number of 
L\"{o}schian distances $l$ in the range $[\,a, b\,]$, the order of the maximum clique $q$, and the check time in the \texttt{kissat} solver (with the \texttt{--forcephase} key that Marijn Heule opened us in secret). 
For most graphs (except for the largest ones), the minimum value of $m$ is given with a step of 5. An asterisk marks the cases for which the dependence of $\chi$ on the position of the bi-chromatic vertex was observed. For $k=9$ and 10, such a study was carried out purposefully, for others it was a by-product of repeated calculations on different solvers. 

In most cases, an obstacle to a further decrease in the value of $d(\chi)$ was the increasing computation time. For $k=7$ and 9, the limiting factor was the growth of $m$ and, as a result, the size 
and the assembly time of CNF files for SAT solver. Suitable values of $m$ were chosen empirically to obtain a $k$-UNSAT solution. On average, the ratio $m/\sqrt{b}$ was about 1.3, but sometimes more than 2 was required. 

\begin{table}[!t]
{
\caption{Exoo type graphs.}
\label{texoo}
\smallskip
\footnotesize

{
\centering

\begin{tabular}{@{}cc@{}}
\begin{tabular}{@{}c@{\;}|@{\;}c@{\;}|@{\;}c@{\;}|@{\:}c@{\:}|@{\:}c@{\:}|@{\;}c@{\;}|@{\;\:}c@{\;\:}|@{\,}r@{}}
\hline
\T\B\
\small{$k$}  & \small{$a$} & \small{$b$} & \small{$d$} & \small{$l$}& \small{$m$} & \small{$q$} & \small{time, s} \\
\hline
\hline \T
 6 &  13 &  21 & 1.27098 &   4 &  5 & 3 & 0.1 \\
   &  19 &  28 & 1.21395 &   5 & 10 & 3 & 1.9 \\
   &  28 &  39 & 1.18019 &   5 & 10 & 3 & 37 \\
   &  49 &  64 & 1.14286 &   6 & 10 & 3 & 1394 \\
   &  73 &  91 & 1.11650 &   7 & 15 & 3 & 41939 \\
   &  91 & 111 & 1.10444 &   8 & 15 & 3 & 169147 \\
   & 169 & 199 & 1.08513 &  11 & 30 & 3 & 4176344 \\
\hline \T
 8 &  27 &  76 & 1.67774 &  18 & 10 & 4 & 0.9 \\
   &  48 & 127 & 1.62660 &  27 & 15 & 4 & 41 \\
   &  57 & 148 & 1.61136 &  30 & 20 & 4 & 89 \\
   &  79 & 201 & 1.59509 &  38 & 20 & 4 & 403 \\
   &  91 & 225 & 1.57243 &  41 & 20 & 4 & 1277 \\
   & 169 & 403 & 1.54422 &  68 & 30 & 4 & 32728 \\
   & 361 & 841 & 1.52632 & 127 & 65 & 4 & 866306 \\
\hline \T
10 &  25 &  97 & 1.96977 &  25 & 15 & 6 &*1.6 \\
   &  36 & 133 & 1.92209 &  33 & 20 & 5 & 47 \\
   &  91 & 331 & 1.90719 &  71 & 20 & 5 & 433 \\
   & 124 & 441 & 1.88586 &  90 & 35 & 5 & 7347 \\
   & 169 & 589 & 1.86687 & 116 & 40 & 5 & 44824 \\
   & 208 & 724 & 1.86568 & 139 & 45 & 5 & 189431 \\
   & 241 & 832 & 1.85803 & 158 & 45 & 5 & 290312 \\
\hline \T
12 &  49 & 247 & 2.24518 &  61 & 20 & 7 & 47 \\
   &  91 & 439 & 2.19640 & 100 & 25 & 7 & 2694 \\
   & 124 & 589 & 2.17945 & 128 & 30 & 7 & 16705 \\
   & 163 & 763 & 2.16356 & 162 & 30 & 7 & 69195 \\
   & 189 & 877 & 2.15412 & 186 & 40 & 7 & 90155 \\
   & 208 & 961 & 2.14946 & 198 & 40 & 7 & 330248 \\
\hline \T
14 &  13 &  93 & 2.67467 &  28 & 10 & 9 & 0.4 \\
   &  25 & 171 & 2.61534 &  47 & 15 & 8 & 833 \\
   &  39 & 256 & 2.56205 &  66 & 20 & 8 & 22768 \\
   &  49 & 304 & 2.49080 &  76 & 30 & 8 & 212177 \\
   &  91 & 549 & 2.45621 & 129 & 30 & 8 & 2263048 \\
\hline
\end{tabular}

&
\begin{tabular}{@{}c@{\;}|@{\;}c@{\;}|@{\,}c@{\,}|@{\:}c@{\:}|@{\:}c@{\:}|@{\;}c@{\;}|@{\;\:}c@{\;\:}|@{}r@{}}
\hline
\T\B\
\small{$k$}  & \small{$a$} & \small{$b$} & \small{$d$} & \small{$l$}& \small{$m$} & \small{$q$} & \small{\,time, s} \\
\hline
\hline \T
 7 & 121 & 268 & 1.48825 &  42 & 25 & 4 & 2.5 \\
   & 183 & 388 & 1.45610 &  59 & 25 & 4 &*7 \\
   & 217 & 444 & 1.43041 &  65 & 30 & 4 & 12 \\
   & 300 & 604 & 1.41892 &  84 & 35 & 3 & 627 \\
   & 507 & 988 & 1.39596 & 125 & 45 & 3 & 6199 \\
   & 567 &1101 & 1.39348 & 139 & 50 & 3 & 7743 \\
   & 675 &1300 & 1.38778 & 159 & 60 & 3 & 15341 \\
\hline \T
 9 &  31 & 111 & 1.89226 &  27 & 25 & 5 &*17 \\
   &  49 & 169 & 1.85714 &  37 & 20 & 5 &*10 \\
   &  91 & 304 & 1.82775 &  63 & 25 & 5 &*40 \\
   & 156 & 516 & 1.81871 & 101 & 40 & 5 &*13 \\
   & 268 & 876 & 1.80794 & 163 & 35 & 5 & 75 \\
   & 361 &1161 & 1.79334 & 209 & 50 & 5 & 100 \\
   & 541 &1731 & 1.78875 & 297 & 55 & 5 & 174 \\
\hline \T
11 &  21 &  91 & 2.08167 &  24 & 15 & 7 & 429 \\
   &  36 & 151 & 2.04803 &  38 & 15 & 7 & 13 \\
   &  49 & 193 & 1.98463 &  45 & 25 & 6 & 23027 \\
   &  63 & 247 & 1.98006 &  57 & 20 & 6 & 12223 \\
   &  73 & 283 & 1.96894 &  63 & 25 & 6 & 125967 \\
   &  91 & 343 & 1.94145 &  75 & 25 & 6 & 307255 \\
   &&&&&&\\
\hline \T
13 &  21 & 133 & 2.51661 &  38 & 15 & 8 & 5 \\
   &  36 & 217 & 2.45515 &  56 & 15 & 8 & 35 \\
   &  49 & 283 & 2.40323 &  70 & 20 & 8 &*1034 \\
   &  91 & 511 & 2.36968 & 119 & 25 & 8 & 2252 \\
   & 144 & 793 & 2.34669 & 176 & 35 & 8 & 105923 \\
   & 169 & 919 & 2.33192 & 201 & 40 & 8 & 143781 \\
\hline \T
15 &  13 &  97 & 2.73158 &  29 & 15 & 9 & 62 \\
   &  28 & 201 & 2.67929 &  55 & 15 & 9 & 850 \\
   &  36 & 252 & 2.64575 &  67 & 20 & 9 & 6769 \\
   &  49 & 337 & 2.62251 &  87 & 25 & 9 & 132533 \\
   &  91 & 624 & 2.61861 & 148 & 30 & 9 & 477925 \\
\hline
\end{tabular}
\end{tabular}

}
}
\end{table}

\section{Divination on coffee grounds}

\paragraph{Extrapolation.}

Now we will delve into the realm of assumptions, predictions and speculations. The parameters of e- and w-graphs will act as coffee grounds. As a modern divination tool, we will use a ruler.

\begin{figure}[H]
\centering
{
\centering
\begin{tabular}{@{}c@{\;}cc@{\;}c@{}}
 6 & \includegraphics[scale=0.35]{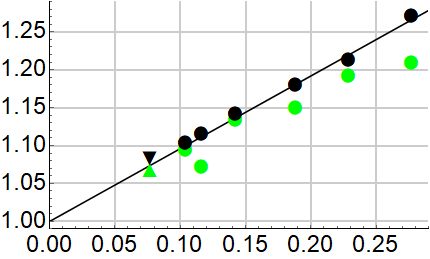}  &  7 & \includegraphics[scale=0.35]{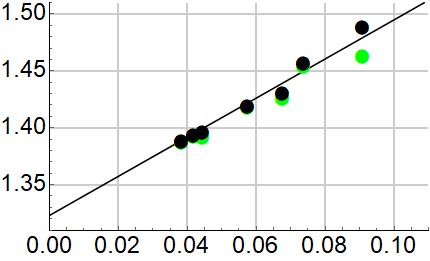}  \\ [1mm]
 8 & \includegraphics[scale=0.35]{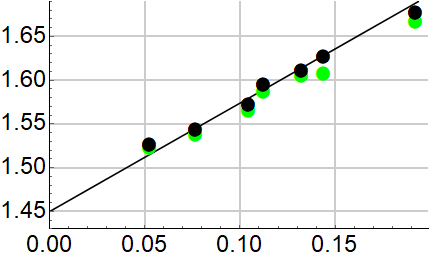}  &  9 & \includegraphics[scale=0.35]{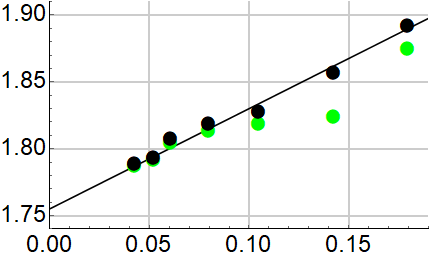}  \\ [1mm]
10 & \includegraphics[scale=0.35]{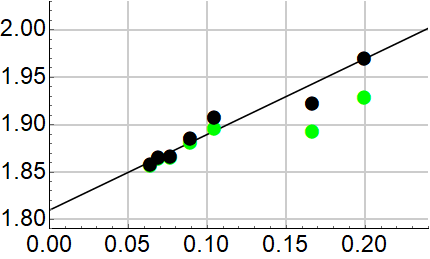} & 11 & \includegraphics[scale=0.35]{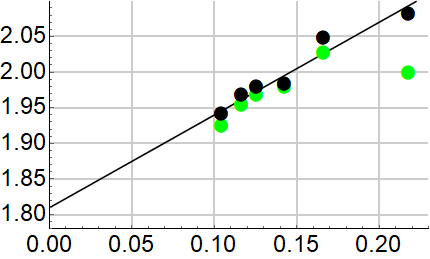} \\ [1mm]
12 & \includegraphics[scale=0.35]{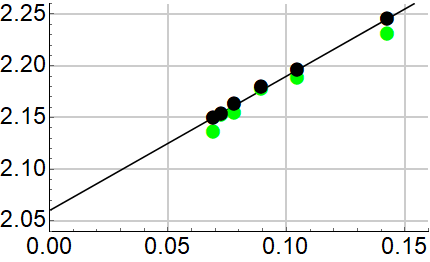} & 13 & \includegraphics[scale=0.35]{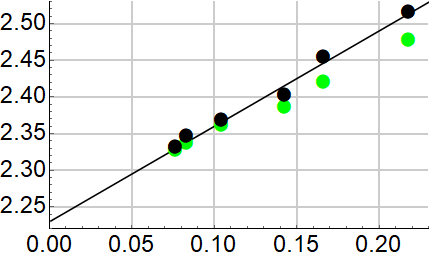} \\ [1mm]
14 & \includegraphics[scale=0.35]{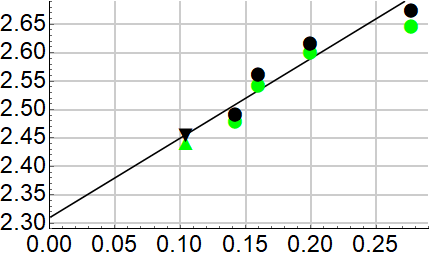} & 15 & \includegraphics[scale=0.35]{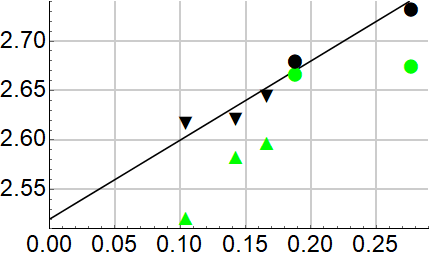} 
\end{tabular} \par
}
\caption{Prediction of asymptotic upper bound of $d$ using estimates $d(r)$ for $6\le k\le 15$. 
Black and green marks correspond to $k$-UNSAT and $k$-SAT solutions for the same $r$.
Inclined lines are drawn along the black marks.}
\label{fpred}
\end{figure}

Let us try to place the estimates $d=\sqrt{b/a}$ given in the Table~\ref{texoo} as points on the coordinate plane depending on the parameter $r=\sqrt{1/a}$, which has the meaning of a resolution element (see Fig.~\ref{fpred}). The dots line up. It can be assumed that with further growth of $a$, new estimates of $d$ 
fall approximately on the same line. This, in turn, allows one to predict the asymptotic value of $d$ (column “pred.” in the Table~\ref{tmain}) for $r\rightarrow 0$. 

In Fig.~\ref{fpred}, $k$-UNSAT solutions are shown in black, and the nearest $k$-SAT solutions for the same $a\in L$ are shown in green. Triangles mark especially hard cases, where there are intermediate values of $b\in L$ for which the solver did not give a solution in an acceptable time.

It is clear that various straight lines can be drawn through a small number of reference points, and the degree of scatter of such predictions will be significant. But the general trend can be caught. Linear extrapolation shows what the upper bound of $d(\chi)$ can supposedly be with unlimited growth of the graph. If we focus on this hypothetical bound, we will not only expand the existing islands of certainty, but also get several new ones. The exceptions in the range $6\le k\le 15$ are 
$k=10$ and 11. 

For w-graphs, we obtained similar predictions (see Table~\ref{tann}) by assuming $r=2\pi/p$. For $k=5$ and 8, the predicted values of $d$ are even smaller than the lower bounds obtained by annulus tiling. This gives an order of 
magnitude of the prediction error. For $k=3$, a w-graph with the smallest possible $d$ has already been found \cite{pmag}.

\paragraph{Confirmations.}
The above scheme of reasoning may seem unreliable, like spring ice. But here are some indirect confirmations that it works: i) the new points do fit quite well on the lines obtained earlier; ii) for small $k$, the predictions practically coincide with the bounds obtained using tilings; iii) for $k=8$, the bound $d=1.444$ was actually predicted a month before we found the corresponding tiling. At that time, the best tiling gave $d=7/5$ \cite{pgray}, and the extrapolation predicted $d\rightarrow 1.44$ (the new points moved this estimate up). Such a good match not only increases the confidence in such predictions, but also allows us to make the following
\begin{conj}
\label{ck8}
The tiling for $k=8$, shown in Fig.~\ref{fper}, is optimal in terms of forbidden distance interval and cannot be improved.
\end{conj}
If we discard caution, then the same can be conjectured for all other $k\le 9$, as well as for $k=12$, meaning the tilings shown in Fig.~\ref{fhex}. By the way, the prediction $d\rightarrow 1$ for $k=6$ can be considered an argument in favor of Exoo's conjecture~\ref{cexoo}. In tilings for $k\ge 13$, we have less confidence.

\paragraph{Inconsistencies.}

As $k$ increases, there are increasing discrepancies between predictions and estimates based on tilings. So, for $k=9$, instead of the expected $d=\sqrt3$ (the existence of a tiling with $d>\sqrt3$ seems incredible), we get $d\rightarrow 1.755$, and for $k=12$, $d\rightarrow 2.05$ is predicted instead of 2. Even more mysterious are the predictions $d\rightarrow 1.81$ for $k=10$ and 11, despite the fact that all our attempts to get a tiling with $d>\sqrt3$ failed. 

Deviations from expectations require explanation. We tend to believe that there are several sources of bias that shift the estimates predicted by the linear extrapolation method towards larger values.
Here are some of them:
i) non-optimal choice of the $q$-clique, as well as the position of the bi-chromatic vertex;
ii) insufficient graph order;
iii) limited possibilities for choosing better graph parameters with increasing $a$ and/or $k$ due to computational difficulties.
Also, we are not sure that the hexagonal lattice used in e-graphs is optimal and leads to unbiased estimates.

Probably, there also exist some as yet unknown constraints, similar to the ($q+3$)-argument using poly-chromatic vertices, which start working at, say, $d>\sqrt3$, reducing prediction bias.

\section[Dancing with a tambourine]{Dancing with a tambourine\footnote{Traditional ritual when the shaman tries to make it rain, but, oddly enough, it doesn't rain.}}

The special status of the cases $k=10$ and 11 (the problem of \textit{redundant} colors, adding of which does not increase $d$), as well as $k=9$ (a noticeable discrepancy between predictions from graphs and tilings), prompted us to conduct additional research. All these cases are connected by the mysterious number $\sqrt3$.

\paragraph{Obscured by clouds.}
Perhaps, we thought, if we take more points, predictions by linear extrapolation will be more reliable and revealing.

\begin{figure}[!t]
\centering
{
\centering
\begin{tabular}{ccc}
    \includegraphics[scale=0.35]{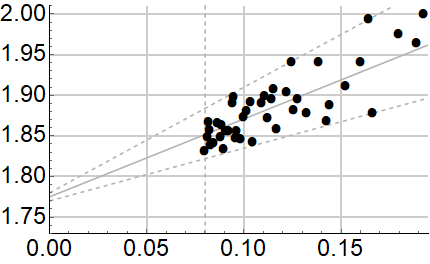} & & \includegraphics[scale=0.35]{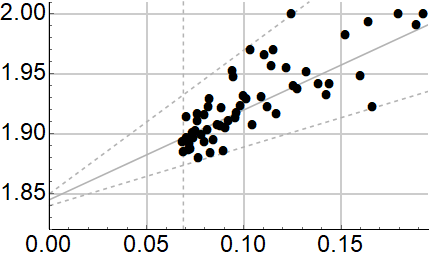} \\ [1mm]
    \includegraphics[scale=0.35]{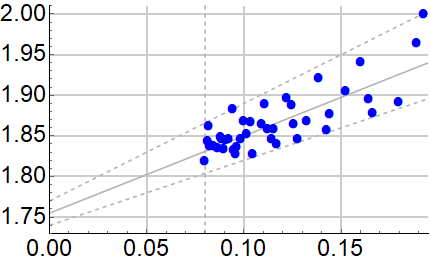} & & \includegraphics[scale=0.35]{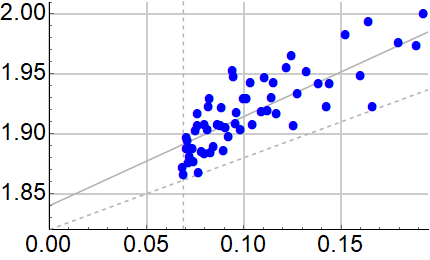} \\ [1mm]
    \includegraphics[scale=0.35]{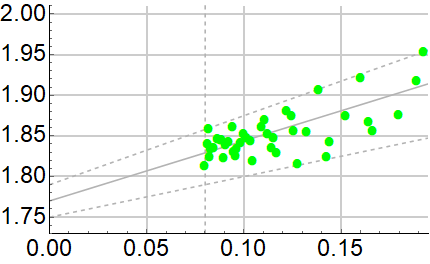} & & \includegraphics[scale=0.35]{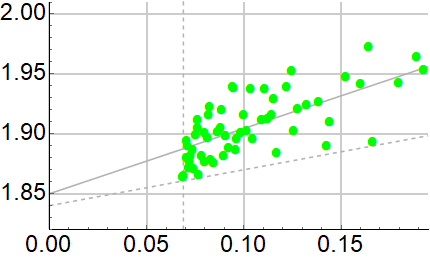} \\
    $k=9$,\;\; $27\le a\le 156$ & & $k=10$,\;\; $27\le a\le 211$ 
\end{tabular} \par
}
\caption{Clouds of e-graph estimates $d(r)$ for $k=9$ and 10. From top to bottom, 
all minimal $k$-UNSAT solutions with a random (black) and optimal (blue) position of the bi-chromatic vertex are shown, as well as $k$-SAT solutions (green) corresponding to the latter case. 
}
\label{fcloud}
\end{figure}

For e-graphs with $k=9$ and 10, we limited the range of $a$ values under study (to reduce the computation time), but for each $a\in L$ in this range, we determined the minimum value of $b$ that led to a $k$-UNSAT solution. As a result, we got a cloud of points (upper graphs in Fig.~\ref{fcloud}). As expected, the dependence $d(r)$ is far from being monotonic. There are both successful pairs $(a, b)$ and pairs that give weak estimates of $d$.

In the previous section, to build a prediction line, we used the optimal values of $r$ (sliding along the bottom edge of the cloud). However, it turned out that if we use all available points to build a straight line (shown in gray in Fig.~\ref{fcloud}), minimizing the standard deviation, the predictions at $r\rightarrow 0$ do not decrease. Perhaps there are not enough points available, or another way to draw a line should be used.

\paragraph{Bi-chromatic vertices.}
Maybe optimizing the position of bi-chromatic vertex will improve the estimates and reduce the prediction bias?

In the same range of $a$ for $k=9$ and 10, we tested all possible distances $s\in[\sqrt{a}, \sqrt{b}]$ between tri- and bi-chromatic vertices. Only some $s$ provide a minimum value of $b$. For example, for $k=9$ we obtain the following optimal sets $(a, b, s^2)$: $(31, 111, 37)$, $(49, 169, 64)$, $(91, 304, 100)$, $(156, 516, 169)$. For some $(a, b)$ there are several such sets. If we do not take into account the position of the bi-chromatic vertex and use an arbitrary $s$ (usually, we took $s=\sqrt{a}$), then this gives weaker estimates: $(31, 121)$, $(49, 171)$, $(91, 309)$, $(156, 523)$.

Thus, in a significant number of cases, the influence of the relative position of the bi-chromatic vertices was confirmed (see Fig.~\ref{fcloud}). But oddly enough, with a general downward shift and a decrease in the spread of $d$ estimates, 
the cloud of points changes so that the predictions at $r\rightarrow 0$ obtained by straight lines drawn along the center and edges of the cloud increase and diverge (do not shrink to a single point, as expected).

Note that the search was not exhaustive: we did not consider all non-isomorphic options of $s$ for placing a bi-chromatic vertex in an e-graph. But this remark concerns only some of $s^2\in L$ 
(for example, 49 and 91).

\begin{figure}[!t]
\centering
{
\centering
\begin{tabular}{@{}c>{\!}c@{}}
    \includegraphics[scale=0.17]{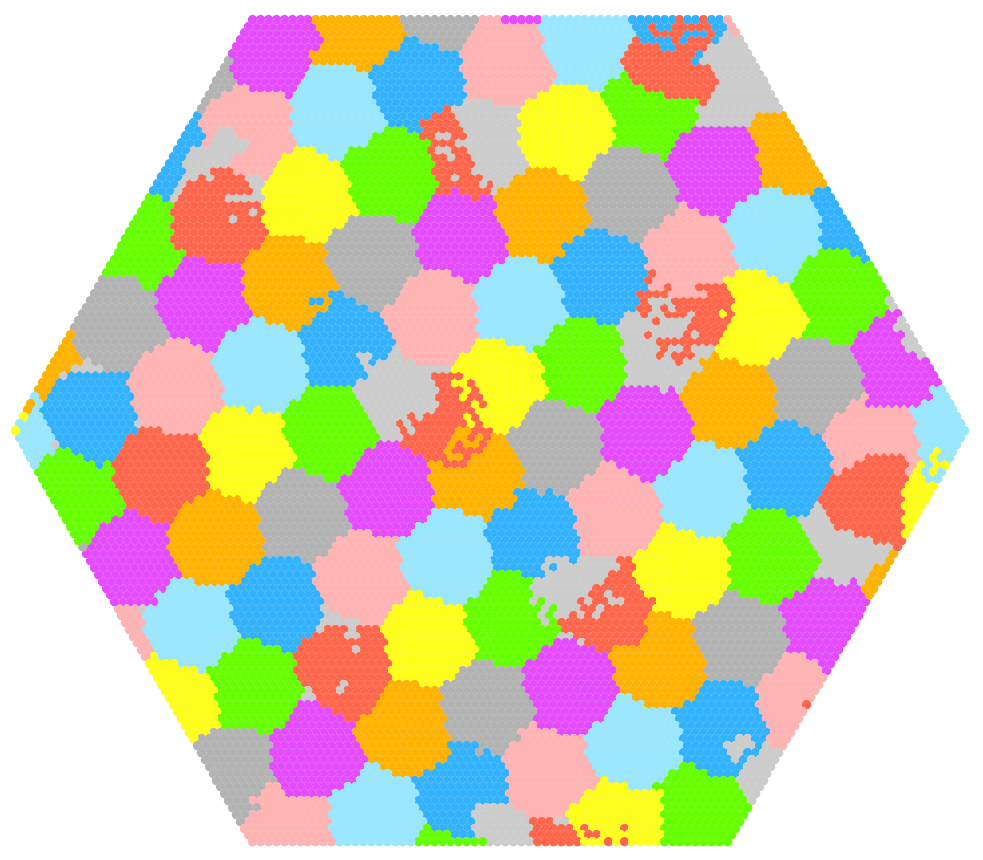} & \includegraphics[scale=0.17]{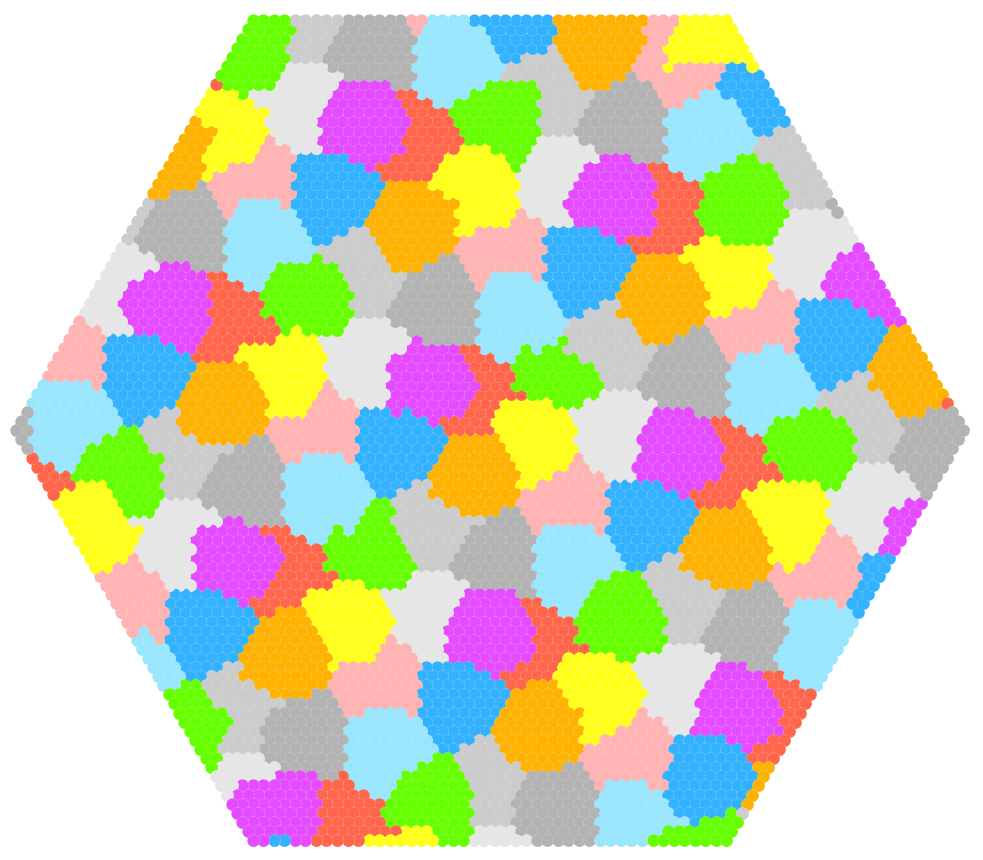} \\
    $k=10$,\; $(m,a,b)=(60,144,508)$ & $k=11$,\; $(m,a,b)=(45,73,279)$ 
\end{tabular} \par
}
\caption{SAT-colorings.}
\label{fsat}
\end{figure}

\paragraph{Impregnable fortress.}
We also tried to go the other way and analyze the obtained $k$-SAT solutions for possible clues in the tiling construction, which would break through the $d=\sqrt3$ threshold for $k=10$ and 11. Fig.~\ref{fsat} shows examples of these SAT colorings. It can be seen that the colors of the vertices are grouped and form pseudo-tiles.

\begin{figure}[!b]
\centering
\includegraphics[scale=0.37]{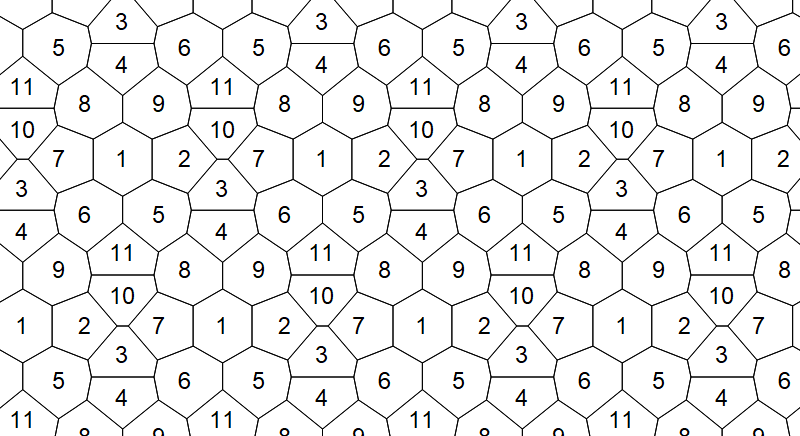}
\caption{An example of an unsuccessful attempt to overcome the magical threshold $d=\sqrt3$ when tiling a plane with 11 colors.
}
\label{fatt}
\end{figure}

But we did not notice significant clues here. In the case of 10 colors, a hexagonal tiling with 9 primary colors and one redundant color is clearly visible. And Fig.~\ref{fatt} shows an adaptation of the considered SAT solution using 11 colors, which also does not allow to beat $d=\sqrt3$. Most of these attempts are broken on the "Mercedes logo": three tiles with a common vertex and unit diagonals, the ends of which form an equilateral triangle with side $\sqrt3$. In Fig.~\ref{fatt}, these are triplets of tiles $\{1, 5, 6\}$ and $\{1, 8, 9\}$.

Let us explain how this constraint works in this case. Denote the common vertices of the triples of tiles $\{1, 5, 6\}$, $\{1, 8, 9\}$, $\{4, 5, 8\}$, $\{4, 6, 9\}$ as $O, A, B, C$. Each of the pairs of tiles $\{1, 5\}$, $\{1, 6\}$, $\{5, 6\}$ in a row is bounded by tiles of the same color, that is, sequences of colors 8158, 9169, 4564 are formed. To obtain the maximum of $d$, points $A, B, C$ must be spaced apart in pairs for the maximum distance. But since the distance from these points to the point $O$ cannot be greater than one, then the smallest of the distances $AB$, $AC$ and $BC$ cannot be greater than $\sqrt3$.
We failed to come up with a tiling for $k=11$ (and even more so for $k=10$) in which the described constraint does not occur, or others do not appear that lead to the same result.

\section{Additions}

\paragraph{Asymptotic estimates}
of $\chi(d)$ for $d\rightarrow\infty$ can be obtained on the basis of the most efficient tiling of the plane and the most dense packing of points.

The densest packing of points on a disk of diameter $d$ with an unlimited increase in the number of vertices $q$ is ensured by their placement at the nodes of a unit hexagonal lattice. The number of points that fit on the disk and give an upper bound of $\chi$ is estimated from the ratio of the areas of the disk and the unit equilateral triangle: $\chi_{ub}\approx\frac{\pi}{\sqrt3}d^2$.

The most efficient tile shape is a regular hexagon with a side of $1/2$, and the best placement of the tile centers is a hexagonal lattice with a step $d+\sqrt3/2$. The required number of colors $\chi$ is obtained by the ratio of the area of the rhombus with side $d$ formed by the lattice to the area of the tile: $\chi_{lb}\approx\frac43d^2$.

Thus, 
$\chi/d^2\in (4/3, \pi/\sqrt3)+o(1)\approx(1.3333, 1.8138)$.
In other words, for large $d$, 
the bounds on $d$ for a fixed $\chi$ differ by a factor of about $7/6$.

\paragraph{Hardware.}
Two computers worked non-stop for a year, checking more than 10\,000 graphs (approximately equal numbers of e- and w-graphs). Some of them required CNF files around 2 GB in size, some took more than a month (unsuccessful with one exception). 

Parameters of computers: i) Intel Core i5-9400F, 2.9 GHz, 6 cores/6 threads, 16 GB RAM; ii) AMD Ryzen5 4600H, 3 GHz, 6 cores/12 threads, 24 GB RAM. Their total performance turned out to be the same. 
Each graph was checked in a separate thread. The computation time is given in terms of the first computer (2 times faster per thread) without taking into account the assembly time of the CNF file.

\paragraph{Observations.}
The \texttt{glucose} and \texttt{kissat} solvers did not differ much in computation time, only on complex graphs requiring many hours, \texttt{kissat} was noticeably more efficient.

In our studies, e-graphs performed better than w-graphs. Apparently, the vertices outside the circle of radius $d$, which we call the \textit{border}, play an important role here. The width of the border (for a fixed $d$) definitely affects both the $k$-colorability and the computational rate, and some optimum is observed.

Only in two cases did the w-graphs prove to be better: i) when estimating $d_{ub}$ for $\chi=9$; ii) when searching for a minimal graph inside the island of stability. In \cite{pmag}, the 18-vertex w-graph allowed us to give a human verifiable proof for $\chi(d)=7$.

The most efficient tiling giving the maximum ratios $d^2/k$ is observed for $k=u^2$, $u\in\mathbb{Z}_{>0}$. In this case, hexagonal tiles of the same color are oriented towards each other with sides (see Fig. 3). A side effect is the inefficiency of subsequent $k$: in the range $k\in[1, 200]$ all $k= u^2+1$ and almost all $k= u^2+2$ give smaller values of $d(k)$, except for $k=6$ and 27. It can be predicted that on all $k= u^2+1$ there will be problems with obtaining tilings that provide $d(k)>d(k-1)$.

\paragraph{Open questions.}
The most obvious (and difficult) goal of subsequent attacks is to break the threshold $d_{lb}=\sqrt3$ for $\chi=10$ and 11, or to prove that this is impossible.
Here are some other questions:

Is it possible to get more efficient tilings by increasing the size of the repeating pattern, for example by using several different tile shapes for each color?
Is it possible to beat the radial coloring of the annulus for $k=6$?
How to explain the frequent repetition of the same values of $a$ and $b$ (for example, 91) in the table of record e-graphs?
What is the main reason for the noticeable difference in the slopes of extrapolation lines for different $k$? (Maybe it's the relative orientation of pseudo-tiles that are assembled from vertices of the same color?)
Is there a bias in the estimates used, and how can it be eliminated?

\paragraph{Conjectures.} 
We have made some progress in confirming Conjecture~\ref{cexoo} by reducing 
$d_{min}$ from 1.285 to 1.085 for $\chi=7$.
Consistent with Conjecture~\ref{cwes}, we have discovered several new islands of certainty, and have optimistic forecasts for several others.
However, we propose the opposite
\begin{conj}
\label{cpar}
For some integers $k\ge 7$, there is no $d$ such that $\chi(d)=k$.
\end{conj}

It means that, for some $\chi\ge 7$, the island of certainty does not exist. In other words, as $d$ grows, the value of $\chi$ can change in non-unit steps. 
Perhaps $k=10$ and 11 are the closest examples of redundant colors for which there is no island of certainty.

\paragraph{Thanks} to Aubrey de Grey for the corrections. Special thanks to Tom Sirgedas for his wonderful program.

In the endless ocean of chromatic numbers it is difficult to find an island of exact knowledge. But if you see several islands at once, this is an occasion to wonder if there is a mainland nearby.


\begin{thebibliography}{100}

\bibitem{wes}
J. Chybowska-Sokół, K. Junosza-Szaniawski, and K. W\k{e}sek.
Coloring distance graphs on the plane. 
\textit{arXiv}:2201.04499v1, 12 Jan 2022.

\bibitem{cou}
D. Coulson. On the chromatic number of plane tilings.
\textit{Journal of the Australian Mathematical Society}, vol.~77, no.~2, 2004, pp.~191--196.

\bibitem{cur}
J. D. Currie and R. B. Eggleton. Chromatic properties of the Euclidean plane.
\textit{arXiv}:1509.03667, 11 Sep 2015.

\bibitem{exoo}
G. Exoo. $\varepsilon$-unit distance graphs.
\textit{Discrete} \& \textit{Computational Geometry}, vol.~33, 
2005, pp.~117--123.

\bibitem{grey}
A.D.N.J. de Grey.
The chromatic number of the plane is at least 5.
\textit{Geombinatorics}, vol.~28, no.~1, 2018, pp.~18--31.

\bibitem{pgray} 
A.D.N.J. de Grey, J. Parts.
Tiling the plane with hexagons: improved separations for $k$-colourings.
\textit{Geombinatorics}, vol.~32, no.~2, 2022, pp.~57--71.

\bibitem{heu}
M.J.H. Heule.
Computing small unit-distance graphs with chromatic number 5.
\textit{Geombinatorics}, vol.~28, no.~1, 2018, pp.~32–-50.

\bibitem{pmag} 
J. Parts.
On the plane and its coloring.
\textit{Geombinatorics}, vol.~31, no.~4, 2022, pp.~189--195.

\bibitem{sir}
T. Sirgedas. 
Hadwiger-Nelson problem / Simulator.
\small \texttt{
https:// groups.google.com/g/hadwiger-nelson-problem/c/a701Kwnhp\_A}
\normalsize

\bibitem{soi}
A. Soifer. \textit{The mathematical coloring book. Mathematics of coloring and the colorful life of its creators}. Springer, New York, 2009. 


\comment{
\bibitem{psma} 
J. Parts.
A small 6-chromatic two-distance graph in the plane.
\textit{Geombinatorics}, vol.~29, no.~3, 2020, pp.~111--115.

\bibitem{plar} 
J. Parts. 
Graph minimization, focusing on the example of 5-chromatic unit-distance graphs in the plane, \textit{Geombinatorics}, vol.~29, No.~4, 2020, pp.~137--166.

\bibitem{pgrn} 
J. Parts. 
What percent of the plane can be properly 5- and 6-colored?
\textit{Geombinatorics}, vol.~30, No.~1, 2020, pp.~25--39.
\textit{arXiv}:2010.12668.
	
\bibitem{ppin} 
J. Parts. 
The chromatic number of the plane is at least 5 – a human-verifiable proof.
\textit{Geombinatorics}, vol.~30, No.~2, 2020, pp.~77--102.
\textit{arXiv}:2010.12661.

\bibitem{pblu} 
J. Parts.
On upper bounds for the multi-fold chromatic numbers of the plane.
\textit{Geombinatorics}, vol.~30, no.~4, 2021, pp.~177--189.

\bibitem{pbro} 
J. Parts.
A 6-chromatic odd-distance graph in the plane.
\textit{Geombinatorics}, vol.~31, no.~3, 2022, pp.~124--137.

\bibitem{pyel} 
J.H. Conway, A.D.N.J. de Grey, J. Parts, A. Soifer. 
Is there a better visualization of Coulson’s "15-colouring of 3-space omitting [monochromatic] distance one"?
\textit{Geombinatorics}, vol.~31, no.~4, 2022, p.~156.

\bibitem{pmag} 
J. Parts.
On the plane and its coloring.
\textit{Geombinatorics}, vol.~31, no.~4, 2022, pp.~189--195.

\bibitem{pgry} 
A.D.N.J. de Grey, J. Parts.
Tiling the plane with hexagons: improved separations for $k$-colourings.
\textit{Geombinatorics}, vol.~32, no.~2, 2022, pp.~57--71.

\bibitem{peme} 
A.D.N.J. de Grey, J. Parts.
On lower bounds of the order of $k$-chromatic unit distance graphs.
\textit{Geombinatorics}, vol.~32, no.~2, 2022, pp.~72--74.



\bibitem{ard}
H. Ardal, J. Ma\v{n}uch, M. Rosenfeld, S. Shelah, and L. Stacho.
The odd-distance plane graph.
\textit{Discrete} \& \textit{Computational Geometry}, vol.~42, no.~2, 2009, pp.~132--141.

\bibitem{cro}
H.T. Croft.
Incidence incidents.
Eureka (Cambridge), No.~30, 1967, pp.~22-26.

\bibitem{exoo}
G. Exoo and D. Ismailescu.
A 6-chromatic two-distance graph in the plane.
\textit{Geombinatorics}, vol.~29, no.~3, 2020, pp.~97--103.

\bibitem{heu2}
M.J.H. Heule.
Odd-distance virtual edges in unit-distance graphs.
\textit{Geombinatorics}, vol.~31, no.~2, 2021, pp.~68--76.

\bibitem{pri}
D. Pritikin.
All unit-distance graphs of order 6197 are 6-colorable.
\textit{Journal of Combinatorial Theory}, Series B 73, 1998, pp.~159--163 

\bibitem{pol}
D.H.J. Polymath.
Comments in Polymath16, thread 17, Dec 2021 -- Feb 2022,
\small \texttt{
https://asone.ai/polymath/index.php?title= Hadwiger-Nelson\_problem
}\normalsize

}

\end{thebibliography}
\end{document}